\documentclass[12pt]{article}
{\begin{Sbox}\begin{minipage}}%
		{\end{minipage}\end{Sbox}\fbox{\TheSbox}}%
\usepackage{color}
\usepackage{amssymb}
\usepackage{amsmath, amsthm, anysize, enumerate, color, url}
\usepackage{graphicx}
\usepackage{subfigure}
\usepackage{xcolor}
\usepackage{caption,subcaption}
\usepackage[colorlinks,citecolor=blue,urlcolor=blue]{hyperref}
\usepackage{booktabs,multirow}
\usepackage{indentfirst}
\newtheorem{theorem}{Theorem} 

\newtheorem{lemma}{Lemma} 
\newtheorem{proposition}{Proposition} 

\newtheorem{definition}{Definition}[section]

\theoremstyle{definition}
\newtheorem{remark}{Remark}  

\usepackage{setspace}
\makeatletter
\renewcommand{\@cite}[1]{[#1]} 
\makeatother

\begin{document}
	\title{ Differential privacy statistical inference for a directed network 
    model with covariates}

\author{
Jing Luo\thanks{Department  of Mathematics and Statistics, South-Central Minzu University, Wuhan, 430074, China, \texttt{Emails:} stajluo@scuec.edu.cn;
$^{\$}$Department of Mathematics and Statistics, South-Central Minzu University, Wuhan, 430074, China, \texttt{Emails:} x1239221809@163.com;$^{\dag}$School of  Statistics and Mathematics, Zhongnan University of Economics and Law ,Wuhan, 430073, China. \texttt{Emails:}qinhong@mail.ccnu.edu.cn; } 
\hspace{15mm}
Hong Qin$^\dag$
\hspace{15mm}
Zhimeng Xu$^{\$}$
\\
 $^{*,\$}$South-Central Minzu University\\ $^{^\dag}$Zhongnan University of Economics and Law
}

	\date{}
	\maketitle

\begin{abstract}
Network data typically contain sensitive relational information, where direct release or sharing may lead to non-negligible privacy violations without proper statistical safeguards.  While differential privacy has emerged as a powerful framework for privacy-preserving network data analysis, theoretical understanding remains limited¡ªparticularly for models incorporating both network structure and nodal attributes. This paper bridges this gap by investigating a directed $\beta$-model with covariates under differential privacy constraints. Our model accounts for both node-level heterogeneity (via $2n$-dimensional degree parameters $\boldsymbol\theta$ ) and covariate-driven homogeneity (via a  $p$-dimensional parameter $\boldsymbol\gamma$). To protect privacy, we introduce a joint Laplace mechanism for releasing network statistics while satisfying differential privacy constraints. Leveraging moment-based estimation techniques, we estimate the parameters of both degree heterogeneity and homogeneity and derive the consistency and asymptotic normality of the differentially private estimators as  the network size tends to infinity. Our theoretical findings are validated through numerical simulations and real-world case studies, demonstrating the validity of our theoretical results.
\end{abstract}


\section{Introduction}
\label{section::introduction}
Both degree heterogeneity and homogeneity are present in the real world. The term ``degree heterogeneity" reflects the network's unequal node degree distribution, where most nodes only keep a small number of edges while a few ``hubs" have extremely high degrees. The significant variations in node degrees are described by this phenomenon. Homogeneity of network vertices refers to the tendency of nodes to establish connections with other nodes that have similar characteristics. The similarity can include the node's attributes, interests, behavioral patterns, etc. For instance, sociologists have found that people who share characteristics (e.g., gender, age, ethnicity, etc.) are more likely to form communities and groups \cite{block2014multidimensional}. This implies that, to some degree, the same covariate information will assist node links in forming. As internet information platforms have grown, a lot of data is being gathered and utilized. Such as healthcare datasets, online social networks, and so forth. Meanwhile, there is also a risk that private information (such as property status, personal preferences, etc.) will be disclosed. Individual rights violations, abuse of personal data, and even legal issues may result from such leaks. Considering that the information about the edges of the network graph and the covariates of the network nodes may involve personal privacy, it is desirable to protect the privacy of the data before performing the relevant analysis.
\par Conventional privacy protection techniques, including data encryption and anonymization, have been shown to be inadequate for protecting private information (e.g., \cite{backstrom2007wherefore}). In order to strictly limit the privacy leakage and to ensure that personal data is not compromised while minimizing the possibility of accessing personal data and increasing its usability, Dwork et al. \cite{dwork2006calibrating} proposed the concept of differential privacy. With this method, adding or deleting records of any individual in a random data release mechanism must not significantly alter the query's output. In the meanwhile, differential privacy offers a mathematically demonstrable privacy-preserving technique that offers a solid theoretical foundation for data exchange and publication.
Up to now, this standard has been extensively adopted as a critical regulatory framework to protect privacy in the course of network data release. (e.g., \cite{lu2017edge,gao2020protecting})
\par While various privacy-preserving algorithms (e.g., \cite{chen2014correlated, karwa2017sharing, salas2023differentially}) exist to secure network data, inferring from noisy data remains highly challenging. For network models, exploring accurate model estimation and asymptotic property analysis of estimators via noisy data is valuable. In recent years, Karwa and Slavkovi{\'c} \cite{karwa2016inference} used denoising methods to study the asymptotic properties of the maximum likelihood estimator (MLE) for the $\beta$-model with discrete Laplace noise and provided detailed proofs. Building on this, Yan \cite{yan2021directed} derived the asymptotic properties of the differential privacy bi-degree sequence maximum likelihood estimation for the $\beta$-model in directed weighted networks, based on both denoised and non-denoised discrete Laplace distributions.
Later, Luo and Qin \cite{luo2022jssc} developed an asymptotic theory for noisy degree sequences. Considering weighted edges in real networks,  Wang et al. \cite{wang2022weighted} demonstrated that the parameter moment estimators constructed based on the bi-degree sequence in the weighted  $p_0$ model possess consistency and asymptotic normality. Unlike Yan, Wang et al.\cite{yan2021directed, wang2022weighted}, who focused their research on one-mode networks, Luo et al. and Wang et al. \cite{luo2022affiliation, wang2022two} investigated two-mode affiliation-weighted networks, deriving the asymptotic properties of differential privacy parameter moment estimators. Concurrently, Luo and Qin \cite{luo2022asymptotic} analyzed the consistency and asymptotic normality of parameter moment estimators for a class of networks with differential privacy degree sequences. Pan and Hu \cite{pan2023differentially} further introduce discrete Laplace noise into the bipartite graph model degree sequence and prove the asymptotic properties of the differential privacy estimator for this model parameter.

\par
Nonetheless, there is still a lack of research on the differential privacy analysis of networks with covariates. Recently, Yan \cite{yan2023differentially} initiated research on node homogeneity under privacy preservation; their focus was limited to undirected networks. As real-world networks, such as trade and conflict networks, are predominantly directed, they are characterized by asymmetric edge relationships and structural complexity (e.g., in- and out-degrees and non-symmetric adjacency matrices). Extending such analyses to directed graphs introduces two key challenges:
(1). Structural asymmetry: Directed edges require modeling bi-directional node interactions, and determining node correspondence and quantifying network similarity are more complex in directed networks than in undirected ones, which increases the difficulty of analysis and theoretical derivation.
(2). Semantic adaptation: Non-symmetric relationships (e.g., information flow direction) demand algorithmic adjustments to preserve both privacy and the directional semantics of network data.
To address this, our paper generalizes the analytical framework of Yan \cite{yan2023differentially} to directed graphs, investigating the statistical theory of covariate-informed parameter estimation for directed network models under differential privacy. We aim to provide rigorous theoretical guarantees for privacy-preserving analyses of real-world directed network data.
\par
We aim to explore the theory of statistical inference for directed $\beta$-models with covariates under differential privacy mechanisms.

Graham \cite{graham2017econometric} first presented the model, and Yan et al. \cite{yan2018statistical} expanded it to its directed version. It comprises a $2n$-dimensional node accessibility parameter $\boldsymbol\theta$ and a p-dimensional regression parameter $\boldsymbol\gamma$ (representing variable effects), where parameter $\boldsymbol\gamma$ reflects the degree of homogeneity (see Section  \ref{subsection 2.1} for details).
The main research contributions are as follows:
First, we employ a jointly Laplace mechanism satisfying differential privacy to release network statistics. As the number of nodes approaches infinity, we derive upper bounds for the error between the differentially private estimators  $(\widehat{\boldsymbol\theta},\widehat{\boldsymbol\gamma})$ and their true values $(\boldsymbol\theta,\boldsymbol\gamma)$. The proof utilizes a two-stage Newton iteration method: first, fixing $\boldsymbol\gamma$ yields upper bounds for the error between $\widehat{\boldsymbol\theta}$ and $\boldsymbol\theta$ ; then leveraging the optimal error bounds of Newton's method under the Newton-Kantorovich condition \cite{KL1948} to derive upper bounds for the error between $\widehat{\boldsymbol\gamma}$ and $\boldsymbol\gamma$. Results show that $\widehat{\boldsymbol\theta}$ converges at a rate of ${\log n/n}^{1/2}$, while $\widehat{\boldsymbol\gamma}$ converges at  $\log n/n$.
Next, an analysis of the asymptotic properties of the differential privacy estimator reveals: the asymptotic variance $\widehat{\boldsymbol\theta}-\boldsymbol\theta^*$ of the heterogeneity parameter $\widehat{\boldsymbol\theta}$ incorporates an additional variance factor, and the asymptotic distribution of the homogeneity parameter $\widehat{\boldsymbol\gamma}$ exhibits a bias term. Finally, these theoretical results are validated through simulation experiments and application to the Lazega dataset \cite{lazega2001collegial}.
\par
The subsequent sections of this paper are organized as follows: Section \ref{section::Model-and-differential-privacy} elaborates on the basic framework of the network model with covariates and the core principles of differential privacy; Section \ref{section::Main-results} presents the key results of this study; Section \ref{section::Numerical-studies} completes the theoretical validation through data examples and simulation studies; Section \ref{section::Discussion} summarizes the research conclusions and provides an outlook for future research. The appendix integrates all theorem statements, proof processes, relevant figures and tables.

\section{Model and differential privacy}
\label{section::Model-and-differential-privacy}
\subsection{Covariate-assisted \texorpdfstring{$\beta$}{beta}-model for directed network}\label{subsection 2.1}

Let ${G_n}$ be a directed graph with $n\ge 2$ nodes, subscripted $"1,\dots,n"$. Let $A=(a_{ij})_{n\times n}$ be the adjacency matrix of ${G_n}$, where $a_{ij}$ is the directed edge from nodes $i$ to $j$ and $a_{ij}\in \{0,1\}$. In this paper, we do not consider the autoloop, so $a_{ii}=0$. Let $d_i=\sum_{j\neq i}^{}a_{ij}$ represent the out-degree of node $i$, and $\mathbf{d}=(d_1,d_2,\dots,d_n)^T$ represent the out-degree sequence of graph $G_n$. Similarly, $b_j=\sum_{i\neq j}^{}a_{ij}$ represent the in-degree of node $j$, and $\mathbf{b}=(b_1,b_2,\dots,b_n)^T$ is the in-degree sequence, the $({d},{b})$ is called a bi-degree sequence. For ease of exposition, define \(\mathbf{g}=(d_1,\dots,d_n,b_1,\dots,b_n)^{T}\).
\par This paper focuses on the degree heterogeneity and homogeneity of network models. Following Graham \cite{graham2017econometric}'s assumption, we set all edges to be mutually independent, with \(a_{ij}\) following a Bernoulli distribution whose probability distribution function is
\begin{equation}\label{1}
	P(a_{ij}=1)={\frac{e^{Z_{ij}^T\boldsymbol{\gamma}+\alpha_i+\beta_j}}{1+{e^{Z_{ij}^T\boldsymbol{\gamma}+\alpha_i+\beta_j}}}},
\end{equation}

Here, parameter \(\alpha_i\) represents the output edge strength of node $i$, while \(\beta_j\) represents the input edge strength of node $j$. The covariate vector $Z_{ij}\in \textrm{R}^{p}$

$Z_{ij}\in \textrm{R}^{p}$ can be constructed based on the similarity or dissimilarity between the p-dimensional attribute vectors \(x_i\) and \(x_j\) of nodes $i$ and $j$. For example, when \(x_i\) represents age, the similarity or dissimilarity between nodes $i$ and $j$ can be measured by \(Z_{ij}=|x_i-x_j|\). The p-dimensional parameter \(\boldsymbol{\gamma}\) serves as the coefficient of the covariate \(Z_{ij}\), reflecting both the connection tendency between nodes $i$ and $j$ and quantifying the degree of homogeneity.

\subsection{Differential privacy}
Given the original database be $D$ (containing $n$ individual data records), $Q$ be a data output randomization mechanism, and \(\mathcal{S}\) represent both the output dataset and the sample space of $Q$. The mechanism \(Q (\cdot|D)\) is a conditional probability distribution over \(\mathcal{S}\). For a positive real number \(\epsilon\), consider two datasets \(D_1\) and \(D_2\) differing by only one record. If the output mechanism \(Q\) satisfies \(\epsilon\)-differential privacy, then for any measurable subset \(S \subseteq \mathcal{S}\), there exists
\[
Q(S |D_1) \leq e^{\epsilon }\times Q(S  |D_2).
\]
The privacy parameter \(\epsilon\) is determined and publicly disclosed by the data administrator to establish the balance between privacy protection and data utility. Smaller $\epsilon$ values mean more privacy. In order to maximize data protection, edge differential privacy is often referenced. Let $\delta(G, G^\prime)$ be the number of edges on which $G$ and $G^\prime$ differ.
The formal definition of edge differential privacy is as follows.
\begin{definition}[Edge differential privacy]
	Let \(\epsilon>0\) be the privacy parameter, and let \(G\) and \(G'\) be any two adjacent graphs differing from each other by a single edge. A certain randomization mechanism
	
	$Q(\cdot |G)$ is $\epsilon$-edge differentially private if
	\[
	\sup_{ G, G^\prime \in \mathcal{G}, \delta(G, G^\prime)=1 } \sup_{ S\in \mathcal{S}}  \log \frac{ Q(S|G) }{ Q(S|G^\prime ) } \le e^\epsilon,
	\]
	here, \(\mathcal{G}\) denotes the collection of all directed graphs with $n$ vertices, while \(\mathcal{S}\) represents the set of all possible outputs.
	
\end{definition}
Dwork et al. \cite{dwork2006calibrating} put forward the definition of global sensitivity as follows.
\begin{definition}[Global Sensitivity]
	Let $f: \mathcal{G} \rightarrow R^k$. The global sensitivity of $f$ is defined as
	\[
	\Delta(f)=\max_{ \delta(G, G^\prime)=1 } \| f(G) - f(G^\prime) \|_1.
	\]
	Where $\|\cdot \|_1$ is the $L_1$ norm.
\end{definition}
Global sensitivity is employed to calculate the maximum \(L_1\) norm distance between two adjacent graphs. The magnitude of noise added in the differential privacy algorithm $Q$ is determined by this sensitivity. It should be noted that when the output consists of network statistics, marginal differential privacy can be achieved by introducing Laplace noise through the Laplace mechanism (as demonstrated in \cite{dwork2006calibrating}), where the added noise is proportional to the global sensitivity of \(f\).

\begin{lemma}
	\label{lemma1}
	Let \(f:\mathcal{G} \to \mathbb{R}^k\), and let \(e_1, \ldots, e_k\) be independent and identically distributed Laplace random variables with  density function \(e^{-|x|/\alpha}/\alpha\). The Laplace mechanism outputs \(f(G)+(e_1, \ldots, e_k)\), which satisfies \(\epsilon\)-edge differential privacy, where \(\epsilon= \Delta(f)/\alpha\).
	
\end{lemma}
When \(f(G)\) takes integer values, the discrete Laplace random variable can serve as the noise, with its probability mass function defined as:

\begin{equation}\label{2}
	{\rm P}(X = x) = \frac{{1 - \alpha }}{{1 + \alpha }}{\alpha ^{\left| x \right|}},x = 0, \pm 1, \ldots ,\alpha  \in (0,1).
\end{equation}
Note that if the discrete Laplace distribution is used, Lemma \ref{lemma1} still holds, and the privacy parameter satisfies \(\epsilon = -\Delta(f)\log \alpha\).
\begin{lemma}[ \cite{yan2021directed}]
	\label{lemma2}
	Let $X$ be a continuous Laplace random variable with the density function $e^{-|x|/\alpha}/\alpha$, then $X$ is a subexponent with parameter $\alpha$; If $X$ be a discrete Laplace random variable with the same probability distribution (\ref{2}), then $X$ is also a subexponent with parameter $2{(\log \frac{1}{\alpha })^{ - 1}}$.
\end{lemma}
\par A nice feature of differential privacy is that any function of the differential privacy mechanism still satisfies differential privacy \cite{dwork2006calibrating}. That is, if $f$ is the output of an $
\epsilon$-differential private mechanism, then $g(f(G))$ is also $\epsilon$-private difference, where $g$ is an arbitrary function.

\subsection{Estimation using differential privacy network statistics}

\par The Laplace mechanism is often used to provide privacy protection in differential privacy problems for network data. By (\ref{1}), we have the  maximum likelihood equation is:
\begin{equation*}
	\begin{aligned}
		&{d}_i=\sum_{j\neq i}^{}{\frac{e^{Z_{ij}^T\boldsymbol{\gamma}+\alpha_i+\beta_j}}{1+{e^{Z_{ij}^T\boldsymbol{\gamma}+\alpha_i+\beta_j}}}},i=1,\dots,n \\
		&{b}_j=\sum_{j\neq i}^{}{\frac{e^{Z_{ij}^T\boldsymbol{\gamma}+\alpha_i+\beta_j}}{1+{e^{Z_{ij}^T\boldsymbol{\gamma}+\alpha_i+\beta_j}}}},j=1,\dots,n-1\\
		&\sum_{j\neq i}^{}Z_{ij}a_{ij}=\sum_{j\neq i}^{}Z_{ij}{\frac{e^{Z_{ij}^T\boldsymbol{\gamma}+\alpha_i+\beta_j}}{1+{e^{Z_{ij}^T\boldsymbol{\gamma}+\alpha_i+\beta_j}}}},
	\end{aligned}
\end{equation*}
where $\mathbf{g} = (\mathbf{b},\mathbf{d}),\mathbf{y}: = \sum\limits_{j \ne i}^{} {{Z_{ij}}} {a_{ij}}$, and $(\mathbf{g},\mathbf{y})$ is the sufficient statistic.
Subsequently, we utilize the continuous Laplace mechanism and its discrete form from Lemma \ref {lemma1} to generate statistics satisfying \(\epsilon\)-edge differential privacy. Given the directed nature of edges between nodes, adding (or removing) a directed edge from node \(i\) to node \(j\) in graph \(G_n\) results in the corresponding graph \(G_n^{'}\).
Then
$${\left\| {\mathbf{g} - \mathbf{g^{'}}} \right\|_1} = 2,{\left\| {\sum\limits_{j \ne i} {({a_{ij}} - a_{ij}^{'}){Z_{ij}}} } \right\|_1} \le pq,$$
In the equation, \(\mathbf{g^{'}}\) denotes the bi-degree sequence of \(G_n^{'}\), and \(a_{ij}^{'}\) represents the weight of \((i,j)\) in \(G_n^{'}\). Letting \(q:=\max_{ijk}|Z_{ijk}|\), we can infer that the global sensitivity of the bi-degree sequence \(\mathbf{g}\) is \({\Delta f_1}=2\), while that of \(\mathbf{y}\) is \({\Delta f_2}=pq\). We present sufficient statistics $\widetilde{\mathbf{g}}=(\widetilde{\mathbf{d}},\widetilde{\mathbf{b}})$ and $\widetilde{\mathbf{y}}$ that satisfy edge-differential privacy as follows:

\begin{equation}\label{3}
	\begin{aligned}
		&\widetilde{d}_i=d_i+ e_i^+,~~i=1, \ldots, n,\\
		&\widetilde{b}_j=b_j+ e_j^{-},~~j=1, \ldots, n,\\
		&\widetilde{y}_t=\sum_{j\neq i}^{}Z_{ijt}a_{ij}+\eta_t,~~t=1,\dots,p.
	\end{aligned}
\end{equation}
where $e_i^+,i=1,\dots,n$ and $e_j^-,j=1,\dots,n$ are independently generated by a discrete Laplace distribution with parameters $\alpha_{n1}=e^{-\epsilon_n/2}$, and $\eta_t,t=1,\dots,p$ are independently generated from the continuous Laplace distribution with parameter $\alpha_{n2}=pq/\epsilon_n$.
\par Because the outer edge of node $i$ to $j$ is also the inner edge of node $j$ to $i$, the sum of the out-degree of the network nodes is equal to the sum of the in-degree. When \((\alpha,\beta)\) is transformed into \((\alpha - c,\beta + c)\), the probability distribution presented in (\ref{1}) remains unchanged.
To identify model parameters, we set $\beta_n=0$ as described in \cite{yan2016asymptotics}.
Let $\boldsymbol\theta=(\alpha_1,\dots,\alpha_n,\beta_{1},\dots,\beta_{n-1})^{T}$, then we use the moment-based estimating equations with noisy degree sequences.
By replacing the degree sequence $ \mathbf{g}=(\mathbf{d},\mathbf{b}) $ and covariate statistic $\mathbf{y}$ with the differential privacy degree sequence $ \widetilde{\mathbf{g}}=(\widetilde{\mathbf{d}},\widetilde{\mathbf{b}}) $ and $\widetilde{\mathbf{y}}$, the estimation equation is:
\begin{equation}\label{4}
	\begin{aligned}
		&\widetilde{d}_i=\sum_{j=1;j\neq i}^{n}{\frac{e^{Z_{ij}^T\boldsymbol{\gamma}+\alpha_i+\beta_j}}{1+{e^{Z_{ij}^T\boldsymbol{\gamma}+\alpha_i+\beta_j}}}},i=1,\dots,n, \\
		&\widetilde{b}_j=\sum_{i=1;i\neq j}^{n}{\frac{e^{Z_{ij}^T\boldsymbol{\gamma}+\alpha_i+\beta_j}}{1+{e^{Z_{ij}^T\boldsymbol{\gamma}+\alpha_i+\beta_j}}}},j=1,\dots,n-1,\\
		&\widetilde{\mathbf{y}}=\sum_{j\neq i}^{}Z_{ij}{\frac{e^{Z_{ij}^T\boldsymbol{\gamma}+\alpha_i+\beta_j}}{1+{e^{Z_{ij}^T\boldsymbol{\gamma}+\alpha_i+\beta_j}}}}.
	\end{aligned}
\end{equation}
The solution to the above system of equations is the differentially private moment estimator, denoted as \((\widehat{\boldsymbol\theta}, {\widehat{\boldsymbol\gamma}})\), where \(\widehat{\boldsymbol\theta}=(\hat{\alpha}_1,\dots,\hat{\alpha}_n,\hat{\beta}_1,\dots,\hat{\beta}_{n-1})^T\) and \(\hat{\beta}_n=0\). A two-step Newton iteration algorithm will be employed to solve the aforementioned system of equations, with the key results provided in Theorem \ref{THEOREM 1}.

\section{Principal  results}
\label{section::Main-results}
Let \(\mathbb{R} = (-\infty, \infty)\) denote the real number field. For a set \(C \subset \mathbb{R}^n\), let \(C^0\) and \(\overline{C}\) denote the interior and closure of \(C\), respectively. For a vector \(\mathbf{x} = (x_1, \dots, x_n)^\top \in \mathbb{R}^n\), use \(\|\mathbf{x}\|_{\infty} = \max_{1 \le i \le n} |x_i|\) to represent the \(\ell_{\infty}\)-norm of \(\mathbf{x}\). Let \(B(\mathbf{x}, \epsilon) = \{\mathbf{y} : \|\mathbf{x} - \mathbf{y}\|_{\infty} \le \epsilon\}\) be a \(\epsilon\)-neighborhood of \(\mathbf{x}\). As for an \(n \times n\) matrix \(J = (J_{i,j})\), \(\|J\|_{\infty}\) denotes the matrix norm that is induced by the \(\|\cdot\|_{\infty}\)-norm defined on vectors in \(\mathbb{R}^n\):
{$$ \|J\|_{\infty}=\max_{\mathbf{x}\neq 0}\frac{\|J\mathbf{x}\|_{\infty}}{\|\mathbf{x}\|_{\infty}}=\max_{1\le i \le n}\sum_{j=1}^{n}| J_{i,j}|.$$}
Denote the matrix norm \(\|\cdot\|\) as \(\|J\| = \max_{i,j} |J_{i,j}|\), and let \({\boldsymbol\theta}^* = (\alpha^*_1, \dots, \alpha^*_n, \beta^*_{1}, \dots, \beta^*_{n-1})^{T}\) represent the true values of parameter \(\boldsymbol\theta\).
Let $\epsilon_{n1}$ and $\epsilon_{n2}$ be two small
positive numbers. When $\boldsymbol\theta\in B(\boldsymbol\theta^*,\epsilon_{n1}), \boldsymbol\gamma\in B(\boldsymbol\gamma^*,\epsilon_{n2})$, we have
$$|Z_{ij}^T\boldsymbol\gamma+\alpha_i +\beta_j|\le \max_{i\neq j}|\alpha^*_i+\beta^*_j|+2\epsilon_{n1}+pq\|{\boldsymbol\gamma}^*\|_{\infty}+pq\epsilon_{n2}:=\rho_n,$$
where $q :=\max_{i,j}\|Z_{ij}\|_{\infty}$, we assume that $q$ is a constant. Consider function $f(x)=\frac{e^x}{(1+e^x)^2}$ decreasing monotonically over $\left[0,\infty \right.\left.\right)$. Define $x=Z_{ij}^T\gamma+\alpha_i+\beta_j$, then On $\left[0,\rho_n\right]$, $f(x)$ decreases monotonically.It is not difficult to verify, when $\boldsymbol\theta \in B(\boldsymbol(\theta)^{*},\epsilon_{n1})$, ${\boldsymbol\gamma}\in B(\boldsymbol\gamma^*,\epsilon_{n2})$,
\begin{equation}\label{5}
	\begin{aligned}
		&m=\frac{{{e^{ {\rho _n}}}}}{{{{\left( {1 + {e^{ {\rho _n}}}} \right)}^2}}} \le \frac{{{e^{Z_{ij}^T\gamma  + {\alpha _i} + {\beta _j}}}}}{{{{\left( {1 + {e^{Z_{ij}^T\gamma + {\alpha _i} + {\beta _j}}}} \right)}^2}}} \le M = \frac{1}{4},\\
		&\left| {\frac{{{e^{Z_{ij}^T\boldsymbol{\gamma} + {\alpha _i} + {\beta _j}}}\left( {1 - {e^{Z_{ij}^T\boldsymbol{\gamma}  + {\alpha _i} + {\beta _j}}}} \right)}}{{{{\left( {1 + {e^{Z_{ij}^T\boldsymbol{\gamma}  + {\alpha _i} + {\beta _j}}}} \right)}^3}}}} \right| \le \frac{{{e^{Z_{ij}^T\boldsymbol{\gamma}  + {\alpha _i} + {\beta _j}}}}}{{{{\left( {1 + {e^{Z_{ij}^T\boldsymbol{\gamma}  + {\alpha _i} + {\beta _j}}}} \right)}^2}}}\left| {\frac{{\left( {1 - {e^{Z_{ij}^T\boldsymbol{\gamma}  + {\alpha _i} + {\beta _j}}}} \right)}}{{\left( {1 + {e^{Z_{ij}^T\boldsymbol{\gamma}  + {\alpha _i} + {\beta _j}}}} \right)}}} \right| \le {M_1} = \frac{1}{4}.
	\end{aligned}
\end{equation}
\subsection{Consistency}
To establish the asymptotic properties of \((\widehat{\boldsymbol\theta},\widehat{\boldsymbol\gamma})\), we introduce a set of functions.\\
\begin{equation}\label{6}
	\begin{aligned}
		F_i(\boldsymbol\theta,\boldsymbol\gamma)&=\widetilde{d}_i-\sum_{j=1,j\neq i}^{n}\frac{e^{Z_{ij}^T\boldsymbol{\gamma}+\alpha_i+\beta_j}}{1+{e^{Z_{ij}^T\boldsymbol{\gamma}+\alpha_i+\beta_j}}},i=1,\dots,n ,\\
		F_{n+j}(\boldsymbol\theta,\boldsymbol\gamma)&=\widetilde{b}_j-\sum_{i=1,i\neq j}^{n-1}\frac{e^{Z_{ij}^T\boldsymbol{\gamma}+\alpha_i+\beta_j}}{1+{e^{Z_{ij}^T\boldsymbol{\gamma}+\alpha_i+\beta_j}}}, j=1,\dots,n-1, \\
		F(\boldsymbol\theta,\boldsymbol\gamma)&=(F_1(\boldsymbol\theta,\boldsymbol\gamma),\dots,F_n(\boldsymbol\theta,\boldsymbol\gamma),F_{n+1}(\boldsymbol\theta,\boldsymbol\gamma),\dots,F_{2n-1}(\boldsymbol\theta,\boldsymbol\gamma))^T.	
	\end{aligned}
\end{equation}
This set of functions is defined based on the differential privacy moment estimation equation. Let \(F_{\boldsymbol\gamma,i}(\boldsymbol\theta)\) denote the value of \(F_i(\boldsymbol\theta,\boldsymbol\gamma)\) when \(\boldsymbol\gamma\) is fixed. Then, \(F_{\boldsymbol\gamma}(\boldsymbol\theta) = (F_{\boldsymbol\gamma,1}(\boldsymbol\theta),\dots,F_{\boldsymbol\gamma,2n-1}(\boldsymbol\theta))^T\).
If the solution to \(F_{\boldsymbol\gamma}(\boldsymbol\theta)=0\) exists, we denote it as \(\widehat{\boldsymbol\theta}_{\boldsymbol\gamma}\). Additionally, the following two functions are defined to study the asymptotic properties of \(\boldsymbol\gamma\):
\begin{equation}\label{7}
	Q(\boldsymbol\theta,\boldsymbol\gamma)=\widetilde{\mathbf{y}}-\sum_{j\neq i}^{}Z_{ij}{\frac{e^{Z_{ij}^T\boldsymbol{\gamma}+\alpha_i+\beta_j}}{1+{e^{Z_{ij}^T\boldsymbol{\gamma}+\alpha_i+\beta_j}}}},
\end{equation}
\begin{equation}\label{8}
	Q_c(\boldsymbol\gamma)=\widetilde{\mathbf{y}}-\sum_{j\neq i}^{}Z_{ij}\frac{e^{Z_{ij}^T\boldsymbol{\gamma}+\widehat{\alpha}_{\gamma,i}+\widehat{\beta}_{\gamma,j}}}{1+{e^{Z_{ij}^T\boldsymbol{\gamma}+\widehat{\beta}_{\gamma,i}+\widehat{\beta}_{\gamma,j}}}}.
\end{equation}
It is easy to see that
$$F(\widehat{\boldsymbol\theta},\widehat{\boldsymbol\gamma})=0,F_{\boldsymbol\gamma}(\widehat{\boldsymbol\theta}_{\boldsymbol\gamma})=0,~Q(\widehat{\boldsymbol\theta},\widehat{\boldsymbol\gamma})=0,~Q_c(\widehat{\boldsymbol\gamma})=0.$$
First, we derive the error bound between $\widehat{\boldsymbol\theta}_{\boldsymbol{\gamma}}$ and ${\boldsymbol\theta}^*$ for $\boldsymbol\gamma\in B({\boldsymbol\gamma}^*,\epsilon_{n2})$. Specifically, this is achieved by constructing the Newton iteration sequence
$\{ {{{\boldsymbol\theta} ^{\left( {k + 1} \right)}}} \}_{k = 0}^\infty$
(with \({\boldsymbol\theta}^*\) as the initial value) and proving that this sequence converges at a geometric rate.
Where ${{\boldsymbol\theta} ^{(k + 1)}} = {{\boldsymbol\theta} ^{(k)}} - {[ {F_{{\boldsymbol\gamma} }^{'}({{\theta} ^{(k)}})}]^{ - 1}}{F_{\boldsymbol\gamma }}({\boldsymbol\theta} ^{(k)})$. The margin of error is described below.
\begin{lemma}\label{lemma3}
	Let $\epsilon_{n1}$ be a positive number and $\epsilon_{n2}=O({(\log n)^{1/2}}/{n^{1/2}}pq)$, if
	$$\tilde{\epsilon}_ne^{6\rho_n}=o(\sqrt{\frac{n}{\log n}}),$$
	then, with probability tending to  one,  $\widehat{\boldsymbol\theta}_{\boldsymbol\gamma}$ exists and satisfies
	$${\left\| {{\widehat{ \boldsymbol\theta }_{\boldsymbol\gamma}} }- {{\boldsymbol\theta} ^*}\right\|_\infty } = O_p\left(e^{3\rho_n}\tilde{\epsilon}_n\sqrt{\frac{\log n}{n}}\right).$$
\end{lemma}
After completing the above symbol definitions, we proceed to conduct an asymptotic analysis of the differential privacy estimator. According to the derivative rule of composite functions,
\begin{equation}\label{9}
	\begin{aligned}
		0 = \frac{{\partial {F_{\boldsymbol\gamma} }({{\widehat {\boldsymbol\theta}} }_{\boldsymbol\gamma }})}{{\partial {\boldsymbol\gamma ^{\rm{T}}}}} = \frac{{\partial F({{\widehat {\theta} }_{\boldsymbol\gamma} },{\boldsymbol\gamma} )}}{{\partial {{\boldsymbol\theta}^{\rm{T}}}}}\frac{{\partial {{\widehat {\boldsymbol\theta} }_{\boldsymbol\gamma} }}}{{\partial {\boldsymbol\gamma} }} + \frac{{\partial F({{\widehat {\boldsymbol\theta} }_{\boldsymbol\gamma} },{\boldsymbol\gamma} )}}{{\partial {{\boldsymbol\gamma}^{\rm{T}}}}} ,     \\
	\end{aligned}
\end{equation}
\begin{equation}\label{10}
	\frac{{\partial {Q_c}({\boldsymbol\gamma} )}}{{\partial {{\boldsymbol\gamma} ^{\rm{T}}}}} = \frac{{\partial Q({{\widehat {\boldsymbol\theta}}_{\boldsymbol\gamma} },{\boldsymbol\gamma} )}}{{\partial {{\boldsymbol\gamma} ^{\rm{T}}}}} = \frac{{\partial Q({{\widehat {\boldsymbol\theta} }_{\boldsymbol\gamma} },{\boldsymbol\gamma} )}}{{\partial {{\boldsymbol\theta} ^{\rm{T}}}}}\frac{{\partial {{\widehat {\boldsymbol\theta} }_{\boldsymbol\gamma }}}}{{\partial {\boldsymbol\gamma} }} + \frac{{\partial Q({{\widehat {\boldsymbol\theta} }_{\boldsymbol\gamma }},{\boldsymbol\gamma} )}}{{\partial {{\boldsymbol\gamma} ^{\rm{T}}}}}.
\end{equation}
The asymptotic behavior of \(\widehat{{\boldsymbol\gamma}}\) is governed by \(Q_{c}^{'}({\boldsymbol\gamma})\). By substituting \(\frac{{\partial {{\widehat {\boldsymbol\theta} }_{\boldsymbol\gamma} }}}{{\partial {\boldsymbol\gamma} }}\) into (\ref{10}) within (\ref{9}), we obtain the expression \({Q_c}({\boldsymbol\gamma} )\) for the partial derivative of ${\boldsymbol\gamma}$:

\begin{equation}\label{11}
	\frac{{\partial {Q_c}({\boldsymbol\gamma} )}}{{\partial {{\boldsymbol\gamma} ^{\rm{T}}}}} = \frac{{\partial Q({{\widehat {\boldsymbol\theta} }_{\boldsymbol\gamma} },{\boldsymbol\gamma} )}}{{\partial {{\boldsymbol\gamma} ^{\rm{T}}}}} - \frac{{\partial Q({{\widehat {\boldsymbol\theta} }_{\boldsymbol\gamma} },{\boldsymbol\gamma} )}}{{\partial {{\boldsymbol\theta} ^{\rm{T}}}}}{\left[ {\frac{{\partial F(\widehat {\boldsymbol\theta} _{\boldsymbol\gamma} ,{\boldsymbol\gamma} )}}{{\partial {{\boldsymbol\theta} ^{\rm{T}}}}}} \right]^{ - 1}}\frac{{\partial F(\widehat {\boldsymbol\theta} _{\boldsymbol\gamma} ,{\boldsymbol\gamma} )}}{{\partial {{\boldsymbol\gamma} ^{\rm{T}}}}}.
\end{equation}
Given that \(\widehat{{\boldsymbol\theta}}_{\boldsymbol\gamma}\) lacks an explicit form, we define the general form of \(\frac{{\partial {Q_c}({\boldsymbol{\gamma}} )}}{{\partial{\boldsymbol{\gamma}} }}\) :

\begin{equation}\label{12}
	H({\boldsymbol\theta} ,{\boldsymbol\gamma}) = \frac{{\partial Q({\boldsymbol\theta} ,{\boldsymbol\gamma} )}}{{\partial {{\boldsymbol\gamma} ^{\rm{T}}}}} - \frac{{\partial Q({\boldsymbol\theta} ,{\boldsymbol\gamma} )}}{{\partial {{\boldsymbol\theta} ^{\rm{T}}}}}{\left[ {\frac{{\partial F({\boldsymbol\theta} ,{\boldsymbol\gamma} )}}{{\partial {{\boldsymbol\theta} ^{\rm{T}}}}}} \right]^{ - 1}}\frac{{\partial F({\boldsymbol\theta} ,{\boldsymbol\gamma} )}}{{\partial {{\boldsymbol\gamma} ^{\rm{T}}}}},
\end{equation}
then, for \({\boldsymbol\theta}\in B({\boldsymbol\theta}^*,\epsilon_{n1})\), the following equation holds:
\begin{equation}\label{13}
	\frac{1}{{{n^2}}}H({\boldsymbol\theta} ,{{\boldsymbol\gamma} ^ * }) = \frac{1}{{{n^2}}}H({{\boldsymbol\theta} ^*},{{\boldsymbol\gamma} ^ * }) + o(1).
\end{equation}
$H(\boldsymbol\theta,\boldsymbol\gamma)$ is the Fisher information matrix of the centralized likelihood function on the $\boldsymbol\gamma$. So we assume $H(\boldsymbol\theta,\boldsymbol\gamma)$ is positively definite. Furthermore, the dimension of \(H(\boldsymbol\theta, \boldsymbol\gamma)\) is fixed, and each element is the sum of \(n(n-1)\) terms. Based on this, we assume there exists a number \({\lambda _n}\) such that:
\begin{equation}\label{14}
	\mathop {\sup }\limits_{{\boldsymbol\theta}  \in B({{\boldsymbol\theta} ^*},{\epsilon _{n1}})} {\left\| {{H^{ - 1}}({\boldsymbol\theta} ,{{\boldsymbol\gamma} ^*})} \right\|_\infty } \le \frac{{{\lambda _n}}}{{{n^2}}}.
\end{equation}
If $n^{-2}H({\boldsymbol\theta} ,{{\boldsymbol\gamma} ^ * })$ converges to a constant matrix, then $\lambda_n$ is bounded.
Then we will show the uniform consistency of $\widehat{{\boldsymbol\theta}}$ and $\widehat{{\boldsymbol\gamma}}$ as follow theorems, whose proofs are given in Appendix.
\begin{theorem}\label{THEOREM 1}
	Suppose  ${\boldsymbol\theta} \in B({\boldsymbol\theta}^*,\epsilon_{n1}), {\boldsymbol\gamma}\in B({\boldsymbol\gamma}^*,\epsilon_{n2})$.  Let $\tilde{\epsilon}_n=1+8/\epsilon_n$ and $\varpi_n=q/\sqrt{\log n}+pq/n\epsilon_n+e^{9\rho_n}(1+8/\epsilon_n)^2$. If
	$$\tilde{\epsilon}_ne^{6\rho_n}=o(\sqrt{\frac{n}{\log n}}),~~~\varpi_n\lambda_n^2e^{9\rho_n}=o(\frac{n}{\log n}).$$
	Then when $n$ tends to infinity, by probability tending to one, the moment estimators $\widehat{\boldsymbol\theta}$,  $\widehat{\boldsymbol\gamma}$ exist and are satisfied
	\begin{align*}
		&{\left\| {{{\widehat{\boldsymbol\theta} } } - {\boldsymbol\theta^*}} \right\|_\infty } =  O_p\left(e^{3\rho_n}\tilde{\epsilon}_n\sqrt{\frac{\log n}{n}}\right),\\
		&{\left\| {\widehat{\boldsymbol\gamma}  - {\boldsymbol\gamma^*}} \right\|_\infty } = {O_p}\left( \frac{{\varpi_n\lambda_n\log n}}{{{n}}}  \right).
	\end{align*}
\end{theorem}
\begin{remark}
	The condition \(\varpi_n\lambda_n^2e^{9\rho_n}=o(\frac{n}{\log n})\) ensures the consistency of the estimator, thereby establishing a crucial trade-off between the privacy parameters \(\epsilon_n\) and \(\rho_n\). Furthermore, under specific conditions, the convergence rate of \(\widehat{\boldsymbol\theta}\) is \((\log n/n)^{1/2}\), while that of \(\widehat{\boldsymbol\gamma}\) is \(\log n/n\).
\end{remark}
\subsection{Asymptotic normality }\label{section3.2}
The core of the asymptotic behavior of the differentially private moment estimator depends on the Jacobian matrix of \(F(\widehat{{\boldsymbol\theta}})\). We observe that this Jacobian matrix has structural characteristics, so we introduce a class of general matrices to characterize its structure. Let $m$ and $M$ be two positive numbers satisfying \(M\ge m\ge 0\); a \((2n-1)\times(2n-1)\) matrix \(V=(v_{i,j})\) is classified as belonging to the class \(L_{n}(m,M)\) if it satisfies the following conditions:
\begin{equation}\label{15}
	\begin{aligned}
		& m \le v_{i,i}-\sum_{j=n+1}^{2n-1}v_{i,j}\le M,i=1,\dots,n-1;v_{n,n}=\sum_{j=n+1}^{2n-1}v_{n,j}\\
		& v_{i,j}=0,~~~~~i,j=1,\dots,n,i\ne j,\\
		& v_{i,j}=0,~~~~~i,j=n+1,\dots,2n-1,i\ne j, \\
		& m\le v_{i,j}=v_{j,i}\le M,~~~~~i=1,\dots,n,j=r+1,\dots,2n-1, \\
		&v_{i,n+i}=v_{n+i,i}=0,~~~~~i=1,\dots,n-1, \\
		&v_{i,i}=\sum_{k=1}^{n}v_{i,k}=\sum_{k=1}^{n}v_{k,i},~~i=n+1,\dots,2n-1 .\\
	\end{aligned}
\end{equation}
Clearly, $V$ has the following structure:
\begin{align*}
	V=
	\begin{pmatrix}
		V_{11} & V_{12}\\
		V_{12}^\top &V_{22}
	\end{pmatrix}.
\end{align*}
Where $V_{11(n\times n)}$ and $V_{22((n-1)\times (n-1))}$ are diagonal matrices, $V_{12}$ is a nonnegative matrix whose elements are positive.
For narration convenience, define $v_{2n,i}=v_{i,2n}:=v_{i,i}-\sum_{j=1,j\ne i}^{2n-1}v_{i,j}$ and $v_{2n,2n}=\sum_{i=1}^{2n-1}v_{2n,i}$. For $i=1,\dots,n-1$, we have $m\le v_{2n,i}\le M$, and $v_{2n,i}=0$ when $i=n,n+1,\dots,2n-1 $, and $v_{2n,2n}=\sum_{i=1}^{2n}v_{2n,i}=\sum_{i=1}^{n}v_{2n,i}$.
Generally, \(V^{-1}\) does not possess a closed-form expression. To estimate \(V^{-1}\), we employ an approximation using the matrix \(S=(s_{i,j})\) (see \cite{yan2016asymptotics}), defined as follows:
\begin{equation}\label{16}
	s_{i,j}=
	\begin{cases}
		\frac{\delta_{i,j}}{v_{i,i}}+\frac{1}{v_{2n,2n}}, &i,j=1,\dots ,n,\\
		-\frac{1}{v_{2n,2n}}, &i=1,\dots,n,j=n+1,\dots ,2n-1,\\
		-\frac{1}{v_{2n,2n}},&i=1,\dots ,2n-1,j=1,\dots ,n,\\
		\frac{\delta_{i,j}}{v_{i,i}}+\frac{1}{v_{2n,2n}},&i,j=n+1,\dots,2n-1.
	\end{cases}
\end{equation}
where $\delta_{i,j}=1$ when $i=j$, otherwise is 0. The upper bound of the approximation error is presented in Proposition \ref{proposition 1}.
\begin{proposition}[\cite{yan2016asymptotics}]\label{proposition 1}
	If $V\in L_{n}(m,M)$ with $M/m = o(n)$, then for large enough $n$,
	$$\|V^{-1}-S\|\le \frac{c_1M^2}{m^3n^2}.$$ \\
	where $c_1$ is a constant that does not depend on $M,m,n$.
\end{proposition}
This paper constructs an asymptotic theory for \(\widehat{\boldsymbol\theta}\) here, with detailed proofs provided in the appendix.
\begin{theorem}\label{THEOREM 2}
	Assume the conditions of Theorem \ref{THEOREM 1} hold and
	$\tilde{\epsilon}_n^2e^{9\rho_n}+\varpi_n^2\lambda_n^2qe^{3\rho_n}\log n = o(n^{1/2}/\log n)$,\\
	(i)  If \(e^{\rho_n}\sqrt{\log n}/\epsilon_n=o(1)\), for any fixed \(k\ge 1\), as \(n \to \infty\), the vector composed of the first $k$ elements of \((\widehat{\boldsymbol\theta}-\boldsymbol\theta^*)\) asymptotically multivariate normal with mean of $0$ and  covariance matrix given by the upper left block matrix \(k\times k\) of \(S^*\). \(S^*\) is the matrix obtained by replacing the parameter \(\boldsymbol\theta\) of matrix $S$ in (\ref{15}) with the true parameter \(\boldsymbol\theta^*\).\\
	(ii) Let $$s_n^2 = Var(\sum\limits_{i = 1}^n {{e_i^+} - \sum\limits_{i = 1}^{n - 1} {{e_{i}^{-}}} } ) = 2(2n - 1)\frac{{{e^{ - \epsilon_n/2}}}}{{{{\left( {1 - {e^{ - \epsilon_n/2 }}} \right)}^2}}}.$$
	If ${s_n}/v_{2n,2n}^{1/2} \to c$ for some constant $c$, then for any fixed $k\ge 1$, as $n \to \infty$, the vector consisting of the first $k$ elements of $(\widehat{\boldsymbol\theta}-\boldsymbol\theta^*)$ is asymptotically $k$-dimensional multivariate normal with mean $0$ and covariance matrix given by
	$$diag(\frac{1}{{{v_{1,1}}}}, \ldots ,\frac{1}{{{v_{k,k}}}}) + (\frac{1}{{{v_{2n,2n}}}} + \frac{{s_n^2}}{{v_{2n,2n}^2}}){{\rm I}_k}{\rm I}_k^T.$$
	Where ${{\rm I}_k}$ is a $k$-dimensional column vector with all elements of $1$.
\end{theorem}
\begin{remark}
	Theorem \ref{THEOREM 2}(i), which ${e^{{\rho _n}}}\sqrt {\log n} /{_n} = o(1)$ holds, requires ${\epsilon_n} \to \infty $. That is, when a small amount of privacy protection is provided by random noise, for any fixed $i$, the asymptotic variance of $\widehat{\boldsymbol\theta}_i$ is the same as the original great likelihood estimator and is $1/{v_{i,i}}$.
\end{remark}
\begin{remark}
	In Theorem \ref{THEOREM 2}(ii), the asymptotic variance of \(\widehat{\boldsymbol\theta}_i\) contains an additional term \(s_n^2/v_{2n,2n}^2\). This arises because its asymptotic expression includes \(\sum\limits_{i = 1}^n {{e_i^+} - \sum\limits_{i = 1}^{n - 1} {{e_{i}^{-}}} }\), whose variance is of order \(ne^{-\epsilon_n/2}\).
	When $\epsilon_n$ becomes small the variance increases rapidly, and the effect on $\widehat{\boldsymbol\theta}_i$ cannot be ignored, which leads to a variance factor.
\end{remark}
Let$(T_{ij})$ be a $(2n-1)$-dimensional vector, where the $i-$th and $(n+j)-$th elements are $1$, and all other elements are $0$.
Define
\begin{align*}
	{S_{{\boldsymbol\theta _{ij}}}}\left( {{\boldsymbol\theta} ,{\boldsymbol\gamma} } \right) = \left( {{a_{ij}} - \frac{{{e^{Z_{ij}^{\rm{T}}{\boldsymbol\gamma}  + {\alpha _i} + {\beta _j}}}}}{{1 + {e^{Z_{ij}^{\rm{T}}{\boldsymbol\gamma}  + {\alpha _i} + {\beta _j}}}}}} \right){T_{ij}}	,~~
	{S_{{\boldsymbol\gamma _{ij}}}}\left( {{\boldsymbol\theta} ,{\boldsymbol\gamma}} \right) = {Z_{ij}}\left( {{a_{ij}} - \frac{{{e^{Z_{ij}^{\rm{T}}\boldsymbol{\gamma}  + {\alpha _i} + {\beta _j}}}}}{{1 + {e^{Z_{ij}^{\rm{T}}\boldsymbol{\gamma}  + {\alpha _i} + {\beta _j}}}}}} \right).
\end{align*}
To simplify notations, write
\begin{align*}
	&V\left( {{\boldsymbol\theta} ,\boldsymbol\gamma } \right) = \frac{{\partial F\left( {{\boldsymbol\theta} ,\boldsymbol\gamma } \right)}}{{\partial {\boldsymbol\theta ^{\rm{T}}}}},~~{V_{Q\boldsymbol\theta }}\left( {{\boldsymbol\theta} ,\boldsymbol\gamma } \right) = \frac{{\partial Q\left( {\boldsymbol\theta,\boldsymbol\gamma } \right)}}{{\partial {\boldsymbol\theta ^{\rm{T}}}}},\\
	&{{\tilde S}_{{\boldsymbol\gamma _{ij}}}}\left( {\boldsymbol\theta ,\boldsymbol\gamma} \right) = {S_{{\boldsymbol\gamma _{ij}}}}\left( {\boldsymbol\theta,\boldsymbol\gamma } \right) - {V_{Q\boldsymbol\theta}}\left( {\boldsymbol\theta ,\boldsymbol\gamma } \right){\left[ {V\left( {\boldsymbol\theta ,\boldsymbol\gamma } \right)} \right]^{ - 1}}{S_{{\boldsymbol\theta _{ij}}}}\left( {\boldsymbol\theta ,\boldsymbol\gamma } \right).
\end{align*}
Define $N=n(n-1)$ and $$\bar H = \mathop {\lim }\limits_{n \to \infty } \frac{1}{N}H\left( {{\boldsymbol\theta ^*},{\boldsymbol\gamma ^*}} \right).$$
Assuming that the limit in the above equation exists, we can obtain the asymptotic distribution of \(\widehat{\boldsymbol\gamma}\). The specific elaboration is as follows:
\begin{theorem}\label{THEOREM 3}
	Assume the condition specified in Theorem \ref{THEOREM 1} is satisfied, if
	$$e^{18\rho_n}q\tilde{\epsilon}_n^4=o(n^{1/2}(\log n)^2),$$
	then
	$$\sqrt N \left( {\widehat {\boldsymbol\gamma}  - {\boldsymbol\gamma ^*}} \right) = {{\bar H}^{ - 1}}{B_*} + {{\bar H}^{ - 1}} \times \frac{1}{{\sqrt N }}\sum\limits_{j \neq i} {{{\tilde S}_{{\boldsymbol\gamma _{ij}}}}\left( {{\boldsymbol\theta ^*},{\boldsymbol\gamma ^*}} \right)}  + {o_p}\left( 1 \right).$$
	When $n\to \infty$, we have $$\sqrt N c^T\left( {\widehat {\boldsymbol\gamma}  - {\boldsymbol\gamma ^*}} \right)\mathop  \to \limits^d N\left( c^T{{{\bar H}^{ - 1}}{B_*},c^T\bar H^{-1}}c \right).$$\\
	Where
	\begin{align*}
		{B_*} = \mathop {\lim }\limits_{n \to \infty } \frac{1}{{2\sqrt N }}\left[ {\sum\limits_{i = 1}^n {\frac{{\sum\limits_{j \ne i} {{Z_{ij}}\frac{{{e^{Z_{ij}^{\rm{T}}{\boldsymbol\gamma ^*} + \alpha _i^* + \beta _j^*}}(1 - {e^{Z_{ij}^{\rm{T}}{\gamma ^*} + \alpha _i^* + \beta _j^*}})}}{{{{\left( {1 + {e^{Z_{ij}^{\rm{T}}{\gamma ^*} + \alpha _i^* + \beta _j^*}}} \right)}^3}}}} }}{{\sum\limits_{j \ne i} {\frac{{{e^{Z_{ij}^{\rm{T}}{\gamma ^*} + \alpha _i^* + \beta _j^*}}}}{{{{\left( {1 + {e^{Z_{ij}^{\rm{T}}{\gamma ^*} + \alpha _i^* + \beta _j^*}}} \right)}^2}}}} }}}  + \sum\limits_{j = 1}^n {\frac{{\sum\limits_{i \ne j} {{Z_{ij}}\frac{{{e^{Z_{ij}^{\rm{T}}{\gamma ^*} + \alpha _i^* + \beta _j^*}}(1 - {e^{Z_{ij}^{\rm{T}}{\gamma ^*} + \alpha _i^* + \beta _j^*}})}}{{{{\left( {1 + {e^{Z_{ij}^{\rm{T}}{\gamma ^*} + \alpha _i^* + \beta _j^*}}} \right)}^3}}}} }}{{\sum\limits_{i \ne j} {\frac{{{e^{Z_{ij}^{\rm{T}}{\gamma ^*} + \alpha _i^* + \beta _j^*}}}}{{{{\left( {1 + {e^{Z_{ij}^{\rm{T}}{\gamma ^*} + \alpha _i^* + \beta _j^*}}} \right)}^2}}}} }}} } \right].
	\end{align*}
\end{theorem}
\begin{remark}
	The limiting distribution of \(\widehat{\boldsymbol\gamma}\) contains a bias term, which needs correction when constructing confidence intervals and conducting hypothesis tests. As noted by Yan et al. (\cite{yan2018statistical}), an analytical bias correction formula applies: \({\widehat {\boldsymbol\gamma} _{bc}} = \widehat {\boldsymbol\gamma} - {H^{ - 1}}(\widehat {\boldsymbol\theta} ,\widehat {\boldsymbol\gamma} )\widehat B/\sqrt {n(n-1)}\), where \(H^{-1}(\widehat {\boldsymbol\theta} ,\widehat {\boldsymbol\gamma} )\) and \(\widehat{B}\) are derived by replacing \(\boldsymbol\gamma\) and \(\boldsymbol\theta\) with their estimators \(\widehat{\boldsymbol\gamma}\) and \(\widehat{\boldsymbol\theta}\) in their expressions.
\end{remark}

\section{Numerical Analyses}
\label{section::Numerical-studies}
  In this section, we utilize numerical simulations to verify the asymptotic results of differential privacy estimators in directed graph networks with covariates, and further illustrate the accuracy of the privacy parameter estimation through the analysis of actual data.
\subsection{Simulation studies}

To verify Theorem \ref{THEOREM 2} and Theorem \ref{THEOREM 3}, we conducted numerical simulations. Following the approach of Yan et al. \cite{yan2018statistical}, parameters were defined  in linear form: setting \(\alpha_{i+1}^* = (n - 1 - i)L/(n - 1)\) for \(i = 0, \ldots,n - 1\), with \(L\) taking values \(0.1\log n,0.2\log n,0.4\log n,\log n\); Set \(\beta_i^* = \alpha_i^*\) (\(i = 1, \ldots, n - 1\)), with \(\beta_n = 0\). The privacy parameter \(\epsilon_n\) simulates three values: \(3\), \(\log n/{n^{1/6}}\), and \(\log n/{n^{1/4}}\).
For node attributes, each node contains a two-dimensional vector \(X_i = (x_{i1}, x_{i2})\): \(x_{i1}\) follows a discrete distribution (taking values 1 or -1 with probabilities 0.3 and 0.7, respectively), while \(x_{i2}\) follows a Beta(2,2) distribution; The edge-layer covariate \(Z_{ij}=(x_{i1}*x_{j1},|x_{i2}-x_{j2}|)^T\) is specified, with the homogeneity parameter set as \(\boldsymbol{\gamma}^*=c(1,1.5)^T\).
Simulations were performed under two scenarios, \(n=100\) and \(n=200\), comprising 6 experimental groups, each repeated 1000 times.
\par According to Theorem \ref{THEOREM 2}, we have ${\hat \xi _{i,j}} = [\hat {{\alpha _i}} - \hat {{\alpha _j}} - (\alpha _i^* - \alpha _j^*)]/{(1/{\hat v_{ii}} + 1/{\hat v_{jj}})^{1/2}}$,${{\hat \zeta }_{ij}}$\\
$ = ({{\hat \alpha }_i} + {{\hat \beta }_j} - \alpha _i^* - \beta _j^*)/{(1/{{\hat v}_{ii}} + 1/{{\hat v}_{n + j,n + j}})^{1/2}}$ and
${\hat \eta _{i,j}} = [\hat {{\beta _i}} - \hat {{\beta _j}} - (\beta _i^* - \beta _j^*)]/{(1/{\hat v_{n + i,n + i}}}$ \\
${+ 1/{\hat v_{n + j,n + j}})^{1/2}}$ follow the standard asymptotic normal distribution, 
where $\hat{v}_{i,i}$ is the estimate of $v_{i,i}$ by replacing $
(\boldsymbol\theta^*,\boldsymbol\gamma^*)$ with $(\widehat{{\boldsymbol\theta}},\widehat{\boldsymbol\gamma})$. For $(i,j)$, chose three pairs of values, respectively $(1,2),(n/2,n+1/2),(n-1,n)$, then we will use the QQ plot and the Kolmogorov-Smirnov (K-S) test to evaluate the asymptotic behavior of the differentially private estimators $\hat\xi_{i,j}$, $\hat \zeta  _{i,j}$ and $\hat{\eta}_{i,j}$, give the $95\%$ confidence interval coverage, the confidence interval length, and the frequency of nonexistence of moment estimates. Since the results are similar, we only give the results of $\hat{\xi}_{ij}$ here. Meanwhile, we also recorded the $95\%$ confidence interval coverage, the length of the confidence intervals, and the mean absolute deviation for $\widehat{\boldsymbol{\gamma}}$. 
 \par
\begin{figure}[tbp]
	\centering
	\subfigure[n=100]{
		\begin{minipage}{1\textwidth}
			\centering
			\includegraphics[width=0.7\textwidth,height=10cm]{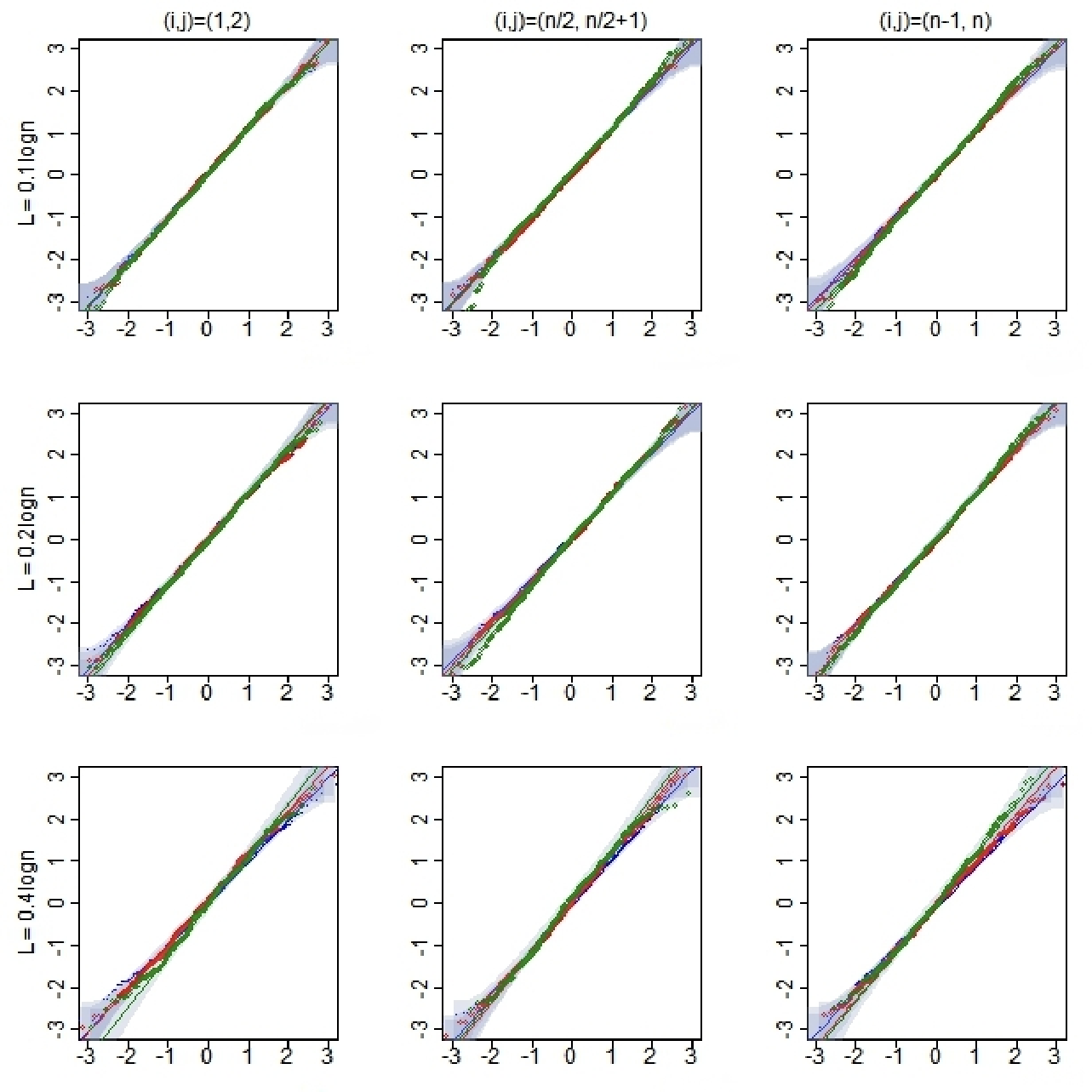}
		\end{minipage}
	}
	\subfigure[n=200]{
		\begin{minipage}{1\textwidth}
			\centering
			\includegraphics[width=0.7\textwidth,height=10cm]{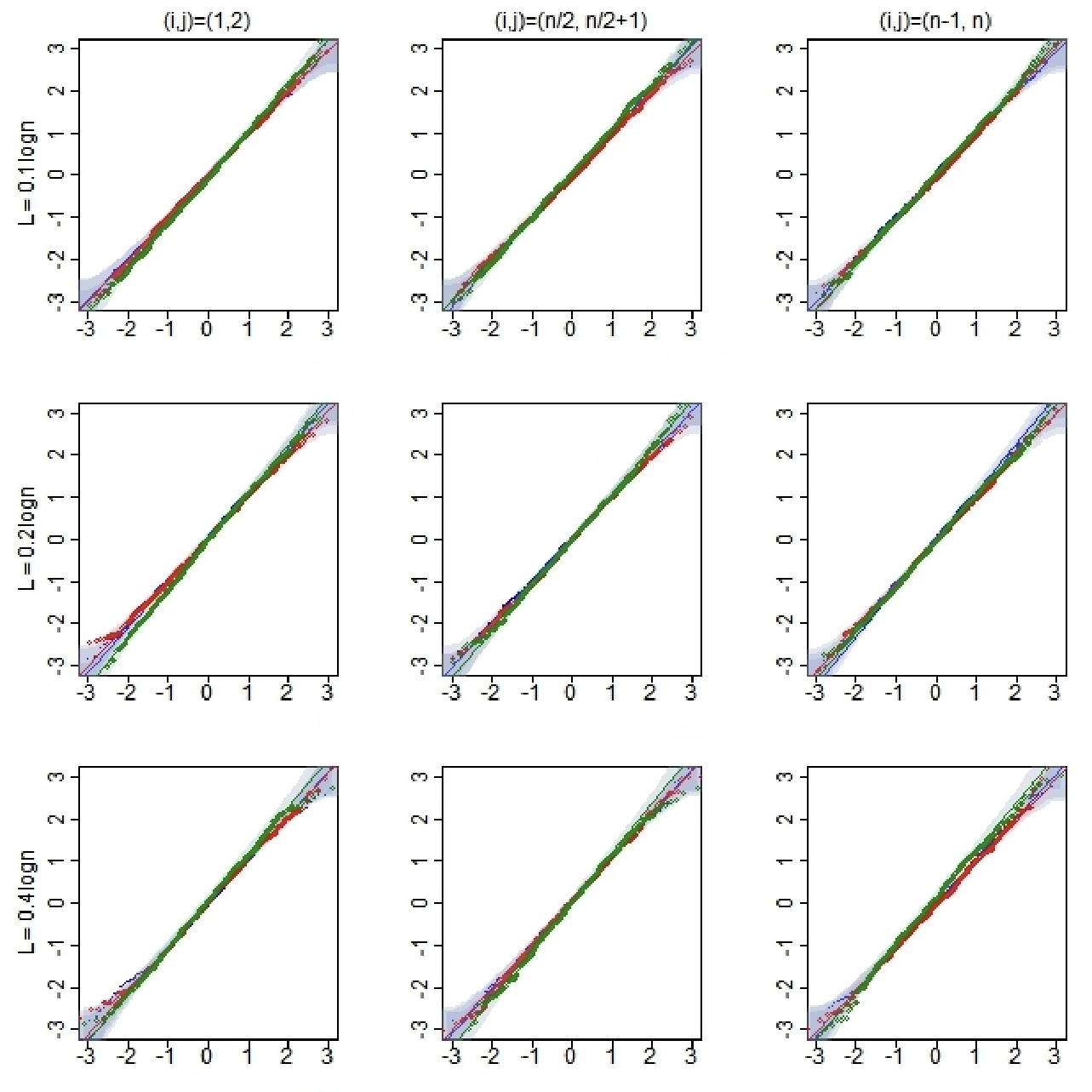}
		\end{minipage}
	}
	\caption{The QQ plot of $\xi_{i,j}$}
	\label{fig-1}
\end{figure}

 
 For the two cases \(n=100\) and \(n=200\), and privacy parameters \(3\), \(\log n/n^{1/6}\), and \(\log n/n^{1/4}\), the simulated QQ-plot is shown in Figure \ref{fig-1}. Different $\epsilon_n$ values are distinguished by blue, red and green. In the figure, the horizontal axis represents the theoretical quantiles, and the vertical axis represents the empirical quantiles. The reference line is $y=x$. Meanwhile, the Kolmogorov-Smirnov (K-S) test statistics and p-values corresponding to each $\epsilon_n$ values are calculated. The results are shown in Table \ref{Table 1}. Combining the QQ-plots and Table\ref{Table 1}: The QQ-plots under the cases of $\epsilon_n=3$ and $\epsilon_n=\log n/n^{1/6}$ have a relatively good fit effect, and as $n$ increases, the overall goodness of fit becomes higher. This indicates that under different combinations of $(i,j)$ and values of $L$, the data distribution has a certain similarity to the theoretical distribution. The K-S test shows that the difference between the data distribution and the normal distribution is not significant, which is consistent with the results of the QQ-plots. Specifically, when \(\epsilon_n=3\), \(L=0.4\log n\), and \((i,j)=({n/2,n/2 + 1})\), the QQ-plot exhibits slight deviations at both ends of the line. In contrast, when \(\epsilon_n=\log n/n^{1/4}\), the deviation of the QQ-plot from the line varies and becomes more pronounced as \(n\) decreases. The p-value from the K-S test confirms this phenomenon. Furthermore, simulations reveal that when \(L=\log n\), the moment estimate does not exist regardless of the value of \(\epsilon_n\), making it impossible to obtain a QQ-plot.

\begin{table}[!ht]
    \centering
    \caption{Kolmogorov-Smirnov (K-S) test statistics and p-values for different privacy parameters  }
	\label{Table 1}
	\small
    \scalebox{0.9}{
	\vspace{0.2cm}
    \begin{tabular}{ccc|cc|cc|cc}
    \hline
       ~& ~&~ &  $\epsilon_n$=3 & ~ & $\epsilon_n=\log n/n^{1/6}$ & ~ & $\epsilon_n=\log n/n^{1/4}$& \\ \hline
       ~& $L$&  $(i,j)$  & K-S & p.value & K-S & p.value & K-S & p.value  \\ \hline
      $n=100$ & $0.1\log n$ & (1,2) & 0.03  & 0.24  & 0.05  & 0.03  & 0.04  & 0.15    \\
     ~& $0.1\log n$ & (50,51) & 0.02  & 0.59  & 0.03  & 0.23  & 0.05  & 0.02    \\
      ~& $0.1\log n$ & (99,100) & 0.02  & 0.89  & 0.02  & 0.82  & 0.03  & 0.49    \\
       ~&  $0.2\log n$& (1,2) & 0.03  & 0.48  & 0.03  & 0.32  & 0.04  & 0.13  \\ 
     ~&$0.2\log n$& (50,51) & 0.03  & 0.24  & 0.03  & 0.21  & 0.03  & 0.33  \\ 
      ~&$0.2\log n$ & (99,100) & 0.03  & 0.34  & 0.03  & 0.30  & 0.03  & 0.37   \\ 
        ~& $0.4\log n$& (1,2) & 0.03  & 0.63  & 0.06  & 0.01  & 0.06  & 0.20    \\ 
     ~&$0.4\log n$& (50,51) & 0.05  & 0.06  & 0.06  & 0.03  & 0.08  & 0.06   \\ 
      ~& $0.4\log n$& (99,100) & 0.06  & 0.34  & 0.03  & 0.30  & 0.03  & 0.37   \\ \hline
      $n=200$ &$0.1\log n$& (1,2) & 0.02  & 0.92  & 0.01  & 1.00  & 0.05  & 0.02   \\ 
      ~&  $0.1\log n$& (100,101)& 0.03  & 0.45  & 0.04  & 0.06  & 0.03  & 0.28   \\ 
       ~& $0.1\log n$ & (199,200) & 0.02  & 0.71  & 0.04  & 0.06  & 0.03  & 0.54   \\ 
      ~& $0.2\log n$ & (1,2) & 0.03  & 0.24  & 0.02  & 0.64  & 0.05  & 0.01   \\ 
       ~& $0.2\log n$ & (100,101)& 0.02  & 0.95  & 0.04  & 0.17  & 0.02  & 0.56   \\ 
        ~& $0.2\log n$ & (199,200) & 0.03  & 0.18  & 0.04  & 0.15  & 0.04  & 0.05   \\
      ~& $0.4\log n$ & (1,2) & 0.04  & 0.15  & 0.03  & 0.35  & 0.05  & 0.09   \\ 
     ~& $0.4\log n$& (100,101) & 0.03  & 0.45  & 0.03  & 0.30  & 0.05  & 0.15   \\ 
      ~& $0.4\log n$ & (199,200) & 0.03  & 0.49  & 0.04  & 0.23  & 0.09  & 0.01   \\ \hline
    \end{tabular}
    }
\end{table}
\begin{table}[h]
	\centering
	\caption{For pair $(i,j)$ of $\alpha_{i}-\alpha_{j}$, the 95\% confidence interval coverage ($\times100\%$), the confidence interval length, and the probabilities of non-existent parameter estimates ($\times100\%$) }
	\label{Table 2}
	\small
	\vspace{0.2cm}
	\begin{tabular}{cccccc}
		\hline
		n         &  $(i,j)$   & $L=0.1\log n$   &  $L=0.2\log n $   &  $L=0.4\log n$& $L=\log n$ \\
		\hline
		\multicolumn {6}{c}{$\epsilon_n$=3}\\
		\bottomrule
		100 & (1,2) & 94.70/1.42/0 & 95.60/1.66/0 & 96.79/2.76/25.20 & NA/NA/100  \\
		~ & (50,51) & 94.20/1.37/0 & 95.60/1.49/0 & 93.18/1.90/25.20 & NA/NA/100  \\ 
		~ & (99,100) & 94.10/1.33/0 & 94.30/1.38/0 & 95.32/1.50/25.20 & NA/NA/100  \\
		200 & (1,2) & 95.20/1.01/0 & 95.20/1.22/0 & 95.84/2.18/8.70 & NA/NA/100  \\
		~ & (100,101) & 94.80/0.97/0 & 95.40/1.07/0 & 94.30/1.42/8.70 & NA/NA/100  \\ 
		~ & (199,200) & 94.70/0.94/0 & 92.80/0.98/0 & 94.19/1.08/8.70 & NA/NA/100  \\ 
		\hline
		\multicolumn {6}{c}{$\epsilon_n=\log n/n^{1/6}$}\\
		\bottomrule
		100 & (1,2) & 93.30/1.42/0 & 94.79/1.66/0.20 & 93.83/3.20/36.80 & NA/NA/100  \\
		~ & (50,51) & 93.00/1.37/0 & 94.49/1.49/0.20 & 92.88/1.90/36.80 & NA/NA/100  \\
		~ & (99,100) & 93.80/1.33/0 & 93.89/1.38/0.20 & 94.78/1.49/36.80 & NA/NA/100  \\
		200 & (1,2) & 94.20/1.01/0 & 95.20/1.22/0 & 94.09/2.18/17.10 & NA/NA/100  \\ 
		~ & (100,101) & 95.80/0.97/0 & 95.00/1.07/0 & 92.88/1.42/17.10 & NA/NA/100  \\ 
		~ & (199,200) & 95.10/0.94/0 & 94.30/0.98/0 & 94.09/1.07/17.10 & NA/NA/100  \\ 
		 \hline
		\multicolumn {6}{c}{$\epsilon_n=\log n/n^{1/4}$}\\
		\bottomrule
	100 & (1,2) & 93.00/1.43/0 & 92.70/1.66/0.50 & 92.55/2.61/71.80 & NA/NA/100  \\
	~ & (50,51) & 92.90/1.38/0 & 92.50/1.49/0.50 & 91.13/1.87/71.80 & NA/NA/100  \\ 
	~ & (99,100) & 91.90/1.34/0 & 92.10/1.39/0.50 & 90.43/1.48/71.80 & NA/NA/100  \\ 
	200 & (1,2) & 92.40/1.01/0 & 92.40/1.22/0 & 91.47/2.19/46.10 & NA/NA/100  \\
	~ & (100,101) & 92.90/0.97/0 & 92.60/1.07/0 & 91.28/1.43/46.10 & NA/NA/100  \\ 
	~ & (199,200) & 93.70/0.94/0 & 94.20/0.98/0 & 92.58/1.08/46.10 & NA/NA/100 \\ \hline
		
	\end{tabular}
\end{table}
\begin{table}[!ht]
	\centering
	\caption{  
    Coverage Probabilities of 95\% confidence interval for \({\gamma}_{i}\) bias-corrected (uncorrected) estimates (x100\%), confidence interval length, average absolute deviation ($\times 10^{-2}$)
    } 
	\label{Table 3}
	\small
	\vspace{0.2cm}
	\begin{tabular}{cccccc}
		\hline
		n         &  $\gamma$   & $L=0.1\log n$   &  $L=0.2\log n $   &  $L=0.4\log n$& $L=\log n$ \\
		\hline
		\multicolumn {6}{c}{$\epsilon_n$=3}\\
		\bottomrule
		 100 &$\gamma_1$ & 94.40(90.70)/0.12/0.89 & 93.00(92.00)/0.14/0.57 & 88.37(90.78)/0.19/0.62 & NA  \\
		~ & $\gamma_2$ & 94.80(93.60)/0.64/5.41 & 94.10(91.30)/0.69/8.08 & 87.70(86.10)/0.85/15.2 & NA  \\ 
		200 &$\gamma_1$& 95.20(91.20)/0.06/0.40 & 90.90(89.60)/0.07/0.17 & 82.69(85.32)/0.10/0.52 & NA  \\
		~ & $\gamma_2$ & 94.20(92.50)/0.32/2.75 & 93.30(89.70)/0.35/4.22 & 84.56(81.27)/0.45/7.85 & NA  \\
		\hline
		\multicolumn {6}{c}{$\epsilon_n=\log n/n^{1/6}$}\\
		\bottomrule
	 100 & $\gamma_1$  & 93.70(89.90)/0.12/0.89 & 92.38(90.48)/0.14/0.57 & 89.56(92.44)/0.19/0.70 & NA  \\ 
	~ &$\gamma_2$  & 93.40(91.90)/0.65/5.42 & 91.18(88.78)/0.69/8.09 & 83.70(83.31)/0.85/16.2 & NA  \\ 
		200 & $\gamma_1$  & 95.60(91.60)/0.06/0.40 & 90.30(89.40)/0.07/0.17 & 81.30(84.32)/0.10/0.51 & NA  \\ 
		~ & $\gamma_2$  & 95.20(93.70)/0.32/2.74 & 93.20(89.90)/0.35/4.21 & 83.35(79.01)/0.45/7.86 & NA  \\ 
		\hline
		\multicolumn {6}{c}{$\epsilon_n=\log n/n^{1/4}$}\\
		\bottomrule
		100 & $\gamma_1$ & 93.00(87.60)/0.13/0.90 & 91.00(88.20)/0.14/0.59 & 85.81(89.72)/0.18/0.55 & NA  \\ 
		~ & $\gamma_2$ & 90.10(88.50)/0.64/5.43 & 88.00(82.10)/0.69/8.10 & 82.27(74.47)/0.85/14.5 & NA  \\ 
		200 & $\gamma_1$& 92.90(89.30)/0.06/0.40 & 90.20(89.50)/0.07/0.17 & 79.96(84.42)/0.10/0.51 & NA  \\ 
		~ & $\gamma_2$ & 92.00(90.90)/0.32/2.75 & 88.30(85.90)/0.35/4.21 & 82.00(77.92)/0.45/7.86 & NA \\ 
		\hline
	\end{tabular}
\end{table}
Table \ref{Table 2} records the 95\% confidence interval coverage, interval length, and differential privacy estimation omission frequency for \(\alpha_i-\alpha_j\). The  findings indicate that: when \(\epsilon_n\) and \(n\) are fixed, the confidence interval length increases as \(L\) increases; when \(\epsilon_n\) and \(L\) are fixed, the interval length decreases as \(n\) increases. Furthermore, when \(\epsilon_n=3\) and \(\log n/n^{1/6}\), the confidence interval coverage rate approaches the 95\% target level; when \(\epsilon_n=\log n/n^{1/4}\), the coverage rate falls slightly below the target level. Note that when \(L=\log n\), differential privacy estimates do not exist with probability one. This shows that setting the parameter space too large can cause errors, so the critical value of \(L\) within this document should be less than \(\log n\).
\par Table \ref{Table 3} records the coverage of the $95\%$ confidence interval, the length of the confidence interval, and the mean of the absolute deviation of the estimator $\widehat{\boldsymbol{\gamma}}$, as well as the deviation correction estimates $\widehat{\boldsymbol{\gamma}}_{bc}(=\widehat{\boldsymbol{\gamma}}-\hat{H}^{-1}\hat{B}/\sqrt{n(n-1)})$. When $\epsilon_n=3,\log n/n^{1/6},L\le 0.2\log n$, the coverage of $\widehat{\boldsymbol{\gamma}}_{bc}$ obtained by bias correction for $\widehat{\boldsymbol{\gamma}}$ is closer to $95\%$ of the target level. When $\epsilon_n=\log n/n^{1/4}$, the confidence interval coverage of $\widehat{\boldsymbol{\gamma}}_{bc}$ is biased. It shows that the selection of the privacy parameter will have some influence on the results of the simulation. In addition, when $n$ is fixed, the length of the confidence interval of $\widehat{\boldsymbol{\gamma}}_{bc}$ increases with the increase of $L$.

\subsection{A data example}
As an example, we use the Lazega dataset of 71 lawyers\cite{lazega2001collegial} , 
which can be downloaded from \href{https://www.stats.ox.ac.uk/~snijders/siena/Lazega\_lawyers\_data.htm}{https://Lazega\_lawyers\_data.htm}. 
This dataset is derived from a web-based study on corporate law partnerships conducted by SG\&R, a corporate law firm located in the northeastern United States, between 1988 and 1991. It comprises network measurements of the firm’s 71 lawyers (including partners and employees), who were surveyed to identify which peers they had sought professional advice in the past year.(The advice was not simply technical advice, but advice that helped them handle their cases correctly and make the right decisions). Besides, the dataset gathers attribute information pertaining to individual lawyers, including formal identity, gender, office location, years of service in the company, age, practice, and law school.
\par
\begin{figure}[h]
	\centering
	\subfigure[]{
		\begin{minipage}[b]{0.3\textwidth}
			\centering
			\includegraphics[scale=0.28]{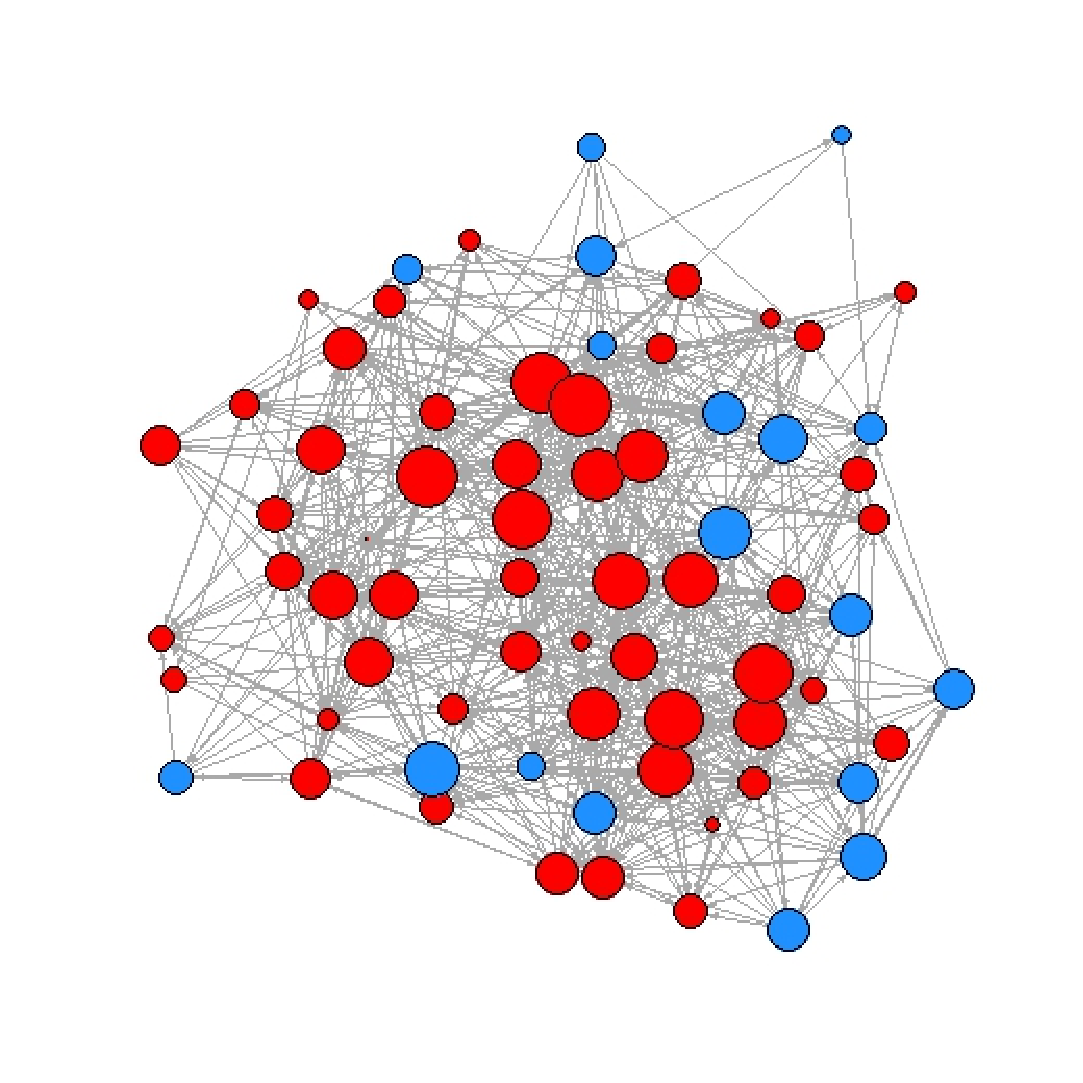}
	\end{minipage}}
	\subfigure[]{
		\begin{minipage}[b]{0.3\textwidth}
			\centering
			\includegraphics[scale=0.3]{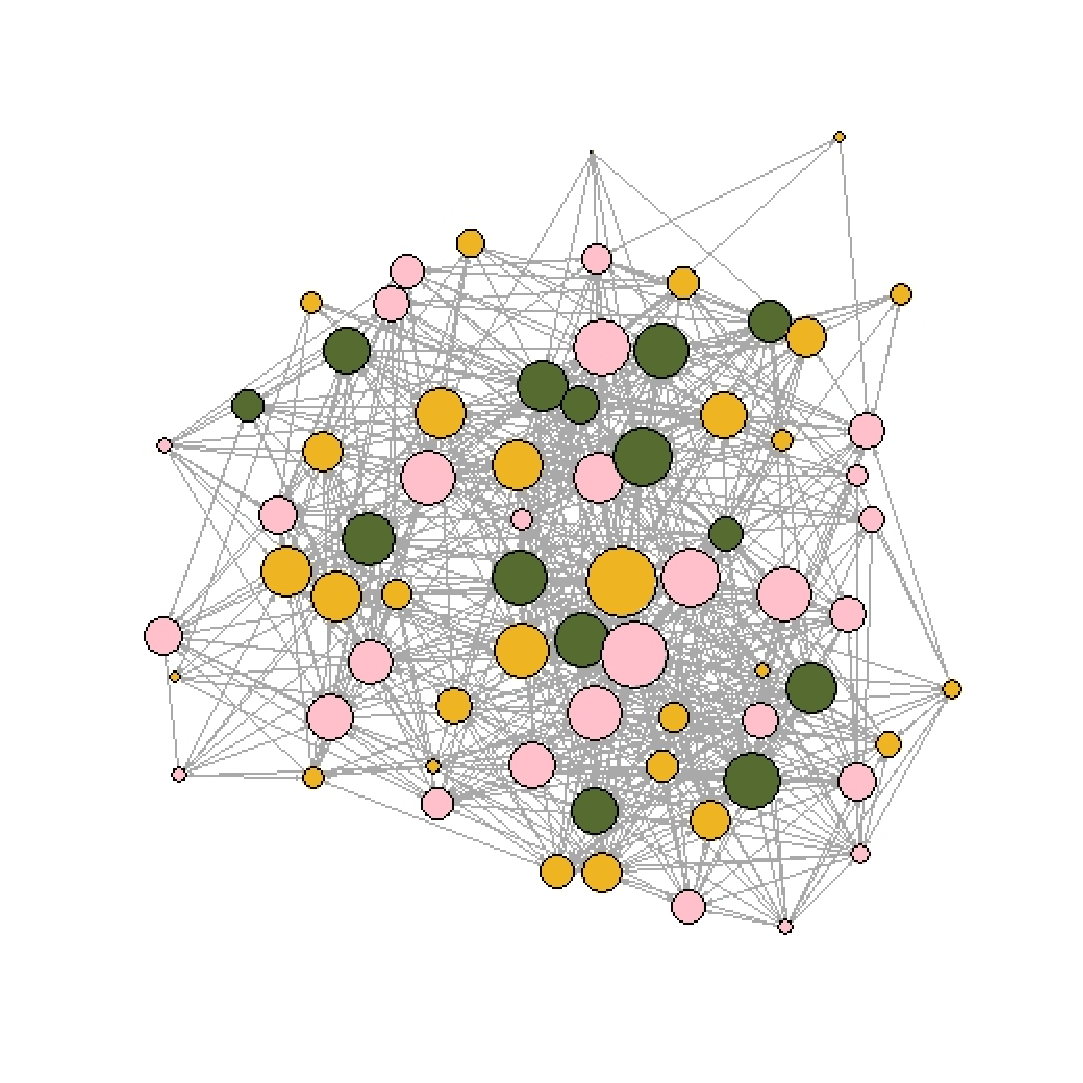}
	\end{minipage}}
	\caption{Visualization of the Lazegas Lawyers’ Advisory Network}
	\label{fig-4}
\end{figure}

 Here, we regard $71$ lawyers as $71$ nodes marked as $1,\dots,71$, and the value of edge $(i,j)$ is $1$ to indicate that lawyer $i$ consulted with lawyer $j$ and received professional advice, otherwise the value of edge $(i,j)$ is $0$. Each node in the dataset has 7 attributes, which are: Formal identity (1=partner; 2=employee); Gender (1=male; 2=female); Office location (1=Boston; 2=Hartford; 3=Providence); Years of service in the company; Age; Practice (1=litigation; 2=company); Law Schools (1: Harvard University, Yale University; 2: University of California; 3: Other). We denote the property of each node $i$ as $X_i=(x_{i1},x_{i2},\dots,x_{i7})$. The covariance \(Z_{ij}\) of edge \((i,j)\) is defined as a function of the attributes of nodes \(i\) and \(j\): for categorical variables, if the \(k\)th attribute of nodes \(i\) and \(j\) is identical, then \(Z_{ijk}=1\); otherwise, \(Z_{ijk}=-1\). For continuous variables, \(Z_{ijk}\) takes the absolute distance \(|x_{ik}-x_{jk}|\). Figure \ref{fig-4}  illustrates a directed network diagram of the dataset containing gender and law school information, where node sizes are proportional to out-degree (Figure (a)) and in-degree (Figure (b)), respectively. In Figure (a), colors distinguish gender (red for males, blue for females), while Figure (b) uses colors to distinguish law schools (green for Harvard and Yale, pink for University of California, orange for other institutions).
\par 
Before analyzing the dataset, we removed node $6$ with an out-degree of $0$ and node $44$ with an in-degree of $0$, proceeding with subsequent analysis on the remaining $69$ nodes. The statistical results for the quartiles, maximum values, and minimum values of the difference quantities $d$ and $b$ are presented in Table \ref {Table 4}.
\begin{table}[h]
	\centering
	\caption{The quartile, maximum and minimum of degree }
	\label{Table 4}
	\small
	\vspace{0.2cm}
	\begin{tabular}{cccccc}
		\bottomrule
		degree& min & $1/4$ quantile & $1/2$quantile & $3/4$ quantile &max\\
		\hline
		d& 2 & 7 & 12 &17 &30\\
		b& 1 &6 &11 &19 &37\\
		\bottomrule
	\end{tabular}
\end{table}\\
\par Then, we substitute the $69\times 69$ adjacency matrix generated by Lazega's advisory network into the model. We set private parameters $\epsilon_n=3$ in this data set, calculate the noisy bi-degree sequence $(\tilde{\mathbf{d}},\tilde{\mathbf{b}})$ of each node, the differentially private estimators of the influencing parameters $\hat{\alpha}$ and $\hat{\beta}$, and the standard error of the estimation, as shown in Table \ref{Table 5}. Where $\hat{\beta}_{71}$ is used as a reference. As can be seen from Table \ref{Table 5}, the heterogeneity parameter estimates of out-degree and in-degree were found to vary widely, ranging from a minimum of $-9.51$ to a maximum of $-2.88$ and a minimum of $-0.41$ to a maximum of $5.08$. The magnitude of these estimates reflects the magnitude of the out-degree of the corresponding node well and can reveal the trend of counseling or being counseled among lawyers. The larger the estimate of the out-degree parameter, the higher the frequency of counseling others; the larger the estimate of the in-degree parameter, the higher the frequency of being counseled by others. 
\par
\begin{figure}[tbp]
	\centering
			\includegraphics[scale=0.5]{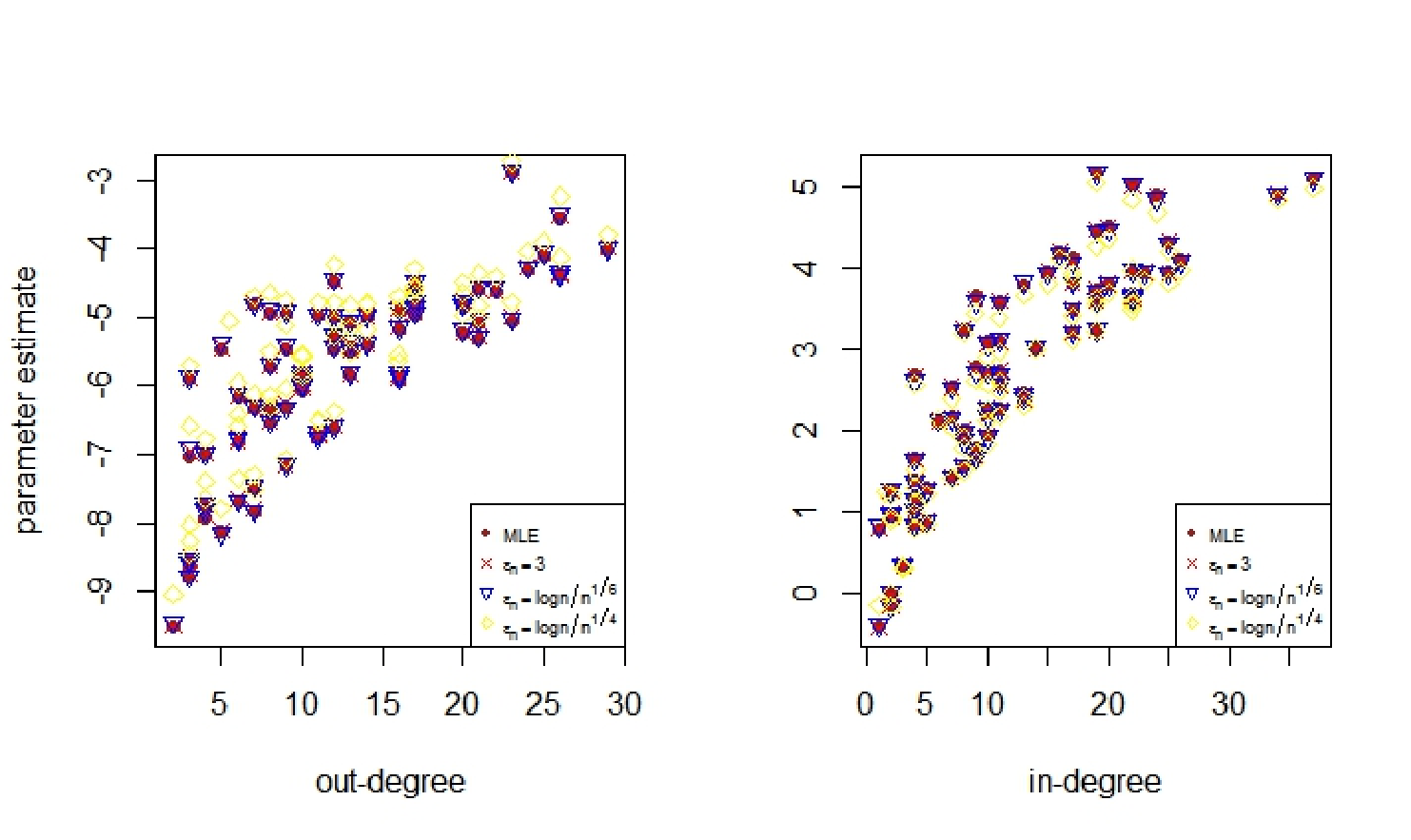}
	\caption{Differential privacy estimation and original maximum likelihood estimation fitting results in the Lazega advisory network dataset}
	\label{fig-5}
\end{figure}

 In addition, we evaluate the difference between the differential privacy estimator $(\widehat{\boldsymbol\theta},\widehat{\boldsymbol\gamma})$ and the original maximum likelihood estimator $(\bar {\boldsymbol\theta} ,\bar{\boldsymbol\gamma } )$. As with the numerical simulations, $\epsilon_n$ was chosen as $3,\log n/{n^{1/6}},\log n/{n^{1/4}}$, respectively. Using the joint Laplace mechanism to repeat the release statistics $g$ and $y$ 1000 times, we computed the values of the differential privacy estimates of the parameters in the Lawyer Advisory Network dataset. For $\epsilon_n=3,\log n/{n^{1/6}},\log n/{n^{1/4}}$, the frequency of non-existence of differential privacy estimates is 49$\%$, 76.9$\%$ and 92.3$\%$, respectively. The specific results are shown in Figure \ref{fig-5} and Table \ref{Table 6}. The vertical axis of the figure represents the estimated value of the parameter $\alpha(\beta)$, and the out (in) degree of the node is shown on the horizontal axis. The black points correspond to $\bar\alpha$ or $\bar{\beta}$ to fit the original data, and the red, blue and yellow points correspond to the median of $\hat{\alpha}$ or $\hat{\beta}$ under different privacy parameters, respectively,  it can be seen from Figure \ref{fig-5}: When $\epsilon_{n}=3,\log n/{n^{1/6}}$,the median value of the differential privacy estimator is very close to the maximum likelihood estimation. When $\epsilon_{n}=\log n/{n^{1/4}}$, the smaller the value of the privacy parameter, the more noise is added. This noise interferes with the estimation of the true parameter, increases the fluctuation of the estimator, and causes the median to deviate. This suggests that the choice of privacy parameters is very important and needs to be considered from a variety of perspectives, an issue that will be pursued subsequently.
\par Table \ref{Table 6} reports the values for the original maximum likelihood  estimate $\bar{\boldsymbol \gamma}$, median values of the differential privacy estimator $\widehat{\boldsymbol\gamma}$, standard errors, confidence intervals and their $p$-values at null hypothesis $\boldsymbol\gamma_i=0(i=1,\dots,7)$. Observations show that: (i). The median value of $\widehat{\boldsymbol\gamma}$ is very close to the maximum likelihood estimation $\bar{\boldsymbol\gamma}$, and lies within the $95\%$ confidence interval;(ii). Five variables—formal identity, gender, office location, years of service at the company, and practice area—exert a significant positive influence on the formation of advisory networks. Specifically, for any of these variables, consistency in their values enhances the likelihood of lawyers engaging in counseling with one another.
Age, on the other hand, negatively affects the formation of counseling relationships among lawyers. In reality, people tend to be more willing to associate with and receive favorable advice from experienced seniors. Furthermore, it is somewhat plausible that having the same law school has little effect on the existence of counseling relationships between lawyers.
\begin{table}[ht]
	\centering
	\caption{The differentially private estimators of $\alpha_i$ and $\beta_j$ and their standard errors in the Lazega's data set when $\hat{\beta}_{71}=0$, private parameter $\epsilon_n=3$.}
	\label{Table 5}
	\small
	\vspace{0.2cm}
	\begin{tabular}{cccccccccccccc}
		\bottomrule
		Vertex & $\tilde{d}_i$ & $\hat{\alpha_i}$ & ${\hat \sigma _i}$ & $\tilde{b}_j$ & $\hat{\beta_i}$ & ${\hat \sigma _j}$ & Vertex & $\tilde{d}_i$ & $\hat{\alpha_i}$ & ${\hat \sigma _i}$ & $\tilde{b}_j$ & $\hat{\beta_i}$ & ${\hat \sigma _j}$  \\ 
		\hline
	1 & 3 & -8.75  & 0.67  & 22 & 3.62  & 0.32  & 37 & 3 & -5.89  & 0.68  & 4 & 2.67  & 0.55   \\ 
	2 & 6 & -7.65  & 0.49  & 22 & 3.96  & 0.32  & 38 & 6 & -6.13  & 0.48  & 17 & 3.49  & 0.38   \\ 
	3 & 6 & -6.75  & 0.53  & 8 & 3.25  & 0.46  & 39 & 14 & -4.98  & 0.36  & 17 & 3.21  & 0.39   \\ 
	4 & 16 & -5.89  & 0.39  & 19 & 3.60  & 0.33  & 40 & 8 & -6.31  & 0.43  & 25 & 3.93  & 0.34   \\ 
	5 & 3 & -8.50  & 0.68  & 17 & 3.81  & 0.37  & 41 & 21 & -4.58  & 0.32  & 22 & 3.63  & 0.35   \\ 
	7 & 4 & -7.73  & 0.63  & 6 & 2.12  & 0.54  & 42 & 25 & -4.07  & 0.31  & 4 & 1.17  & 0.55   \\ 
	8 & 2 & -9.51  & 0.78  & 13 & 2.39  & 0.36  & 43 & 14 & -4.98  & 0.36  & 10 & 2.27  & 0.42   \\ 
	9 & 3 & -8.63  & 0.64  & 14 & 3.00  & 0.37  & 45 & 10 & -5.82  & 0.38  & 8 & 2.00  & 0.43   \\ 
	10 & 7 & -7.50  & 0.47  & 11 & 2.56  & 0.41  & 46 & 7 & -4.80  & 0.46  & 9 & 3.62  & 0.45   \\ 
	11 & 5 & -8.12  & 0.55  & 19 & 3.23  & 0.33  & 47 & 5 & -5.44  & 0.50  & 1 & 0.80  & 1.03   \\ 
	12 & 20 & -5.22  & 0.38  & 20 & 3.80  & 0.33  & 48 & 17 & -4.49  & 0.34  & 4 & 1.38  & 0.56   \\ 
	13 & 16 & -5.82  & 0.37  & 34 & 4.91  & 0.30  & 49 & 9 & -6.28  & 0.42  & 9 & 1.69  & 0.43   \\ 
	14 & 13 & -5.82  & 0.41  & 16 & 4.19  & 0.40  & 50 & 9 & -4.91  & 0.42  & 11 & 3.58  & 0.41   \\ 
	15 & 9 & -5.43  & 0.43  & 19 & 5.14  & 0.34  & 51 & 23 & -2.88  & 0.32  & 2 & 1.25  & 0.79   \\ 
	16 & 26 & -4.38  & 0.36  & 19 & 3.73  & 0.33  & 52 & 13 & -5.48  & 0.37  & 13 & 2.46  & 0.40   \\ 
	17 & 21 & -5.06  & 0.37  & 25 & 4.31  & 0.31  & 53 & 4 & -6.98  & 0.54  & 4 & 0.85  & 0.57   \\ 
	18 & 4 & -7.89  & 0.62  & 11 & 3.09  & 0.42  & 54 & 14 & -5.38  & 0.36  & 9 & 1.78  & 0.44   \\ 
	19 & 29 & -3.98  & 0.35  & 11 & 2.73  & 0.40  & 55 & 22 & -4.60  & 0.32  & 8 & 1.57  & 0.45   \\ 
	20 & 11 & -6.72  & 0.41  & 22 & 3.67  & 0.32  & 56 & 20 & -4.84  & 0.33  & 10 & 1.94  & 0.41   \\ 
	21 & 11 & -6.70  & 0.41  & 22 & 3.62  & 0.32  & 57 & 12 & -5.27  & 0.38  & 11 & 2.26  & 0.42   \\ 
	22 & 12 & -6.59  & 0.42  & 23 & 3.95  & 0.32  & 58 & 11 & -4.97  & 0.39  & 4 & 1.66  & 0.59   \\ 
	23 & 9 & -7.15  & 0.44  & 8 & 1.90  & 0.42  & 59 & 8 & -4.92  & 0.43  & 2 & 0.95  & 0.91   \\ 
	24 & 21 & -5.28  & 0.35  & 26 & 4.08  & 0.31  & 60 & 12 & -4.99  & 0.37  & 7 & 2.16  & 0.45   \\ 
	25 & 8 & -6.51  & 0.47  & 10 & 3.07  & 0.45  & 61 & 3 & -6.97  & 0.63  & 1 & -0.41  & 1.05   \\ 
	26 & 23 & -5.00  & 0.34  & 37 & 5.08  & 0.31  & 62 & 10 & -5.81  & 0.38  & 4 & 1.04  & 0.57   \\ 
	27 & 20 & -4.76  & 0.35  & 9 & 2.76  & 0.39  & 63 & 12 & -4.47  & 0.39  & 2 & 0.98  & 0.79   \\ 
	28 & 26 & -3.51  & 0.39  & 22 & 5.01  & 0.36  & 64 & 8 & -5.70  & 0.42  & 10 & 2.72  & 0.40   \\ 
	29 & 13 & -5.41  & 0.41  & 17 & 4.09  & 0.35  & 65 & 24 & -4.27  & 0.31  & 7 & 1.44  & 0.47   \\ 
	30 & 16 & -5.14  & 0.40  & 20 & 4.50  & 0.34  & 66 & 26 & -4.35  & 0.31  & 2 & -0.15  & 0.75   \\ 
	31 & 17 & -4.84  & 0.39  & 15 & 3.95  & 0.37  & 67 & 16 & -4.88  & 0.35  & 3 & 0.35  & 0.63   \\ 
	32 & 10 & -6.04  & 0.45  & 19 & 4.46  & 0.34  & 68 & 10 & -5.98  & 0.41  & 5 & 0.88  & 0.53   \\ 
	33 & 17 & -4.94  & 0.40  & 7 & 2.53  & 0.48  & 69 & 12 & -5.43  & 0.38  & 3 & 0.34  & 0.63   \\ 
	34 & 6 & -6.79  & 0.51  & 24 & 4.88  & 0.32  & 70 & 7 & -6.30  & 0.43  & 5 & 1.29  & 0.52   \\ 
	35 & 17 & -4.81  & 0.39  & 13 & 3.77  & 0.41  & 71 & 13 & -5.07  & 0.37  & 2 & 0.00  & 0.75   \\ 
	36 & 7 & -7.80  & 0.50  & 10 & 2.19  & 0.40  & ~ & ~ & ~ & ~ & ~ & ~ &   \\ 
		\bottomrule
	\end{tabular}
\end{table}\\
\begin{table}[!ht]
	\centering
	\caption{The original maximum likelihood estimates, differential privacy estimates (Median), standard errors, confidence intervals, and P-Values under the null hypothesis $\gamma_i=0(i=1,\dots,7)$ for $\gamma_i$ in the Lazega advisory network data $(\epsilon_n=3)$}
 
	\label{Table 6}
	\small
	\vspace{0.2cm}
	\begin{tabular}{cccccc}
		\bottomrule
		Covariate & $\bar{\gamma_i}$ & $\hat{\gamma_i}$ &$95\%$ confidence interval&  ${\hat \sigma _i}$ & p-value  \\ 
		\hline
		  identity & 1.137 & 1.124 & [0.975,1.274] & 0.076 & $<$0.001  \\ 
		gender & 	0.523 & 0.520 & [0.389,0.651] & 0.067 & $<$0.001  \\ 
	location & 	1.116 & 1.113 & [0.980,1.246] & 0.068 & $<$0.001  \\ 
		years &	0.055 & 0.054 & [0.032,0.077] & 0.012 & $<$0.001  \\ 
		age & 	-0.029 & -0.029 & [-0.048,-0.011] & 0.009 & $<$0.001  \\ 
			practice &0.953 & 0.949 & [0.843,1.057] & 0.055 & $<$0.001  \\ 
			school & 0.026 & 0.027 & [-0.075,0.129] & 0.052 & 0.465 \\ 
		\bottomrule
	\end{tabular}
\end{table}

\section{Discussion}
\label{section::Discussion}
The present study investigates the asymptotic properties of the heterogeneous parameter $\widehat{\boldsymbol\theta}$ and the homogeneity parameter $\widehat{\boldsymbol\gamma}$ in the directed $\beta$-model with covariates under the differential privacy mechanism. We establish the consistency of differential privacy moment estimates by releasing sufficient network statistics through the joint Laplace mechanism. Furthermore, we derive the asymptotic properties of differentially private estimators $\widehat{\boldsymbol\theta}$ and $\widehat{\boldsymbol\gamma}$ utilizing two-stage Newton iterations. Taking into account the parameter $\widehat{\boldsymbol\theta}$, it is evident that in the presence of noise variance, an additional variance factor is introduced into the asymptotic variance of $\widehat{\boldsymbol\theta}-\boldsymbol\theta^*$
. The results of numerical simulations demonstrate that the theoretical results can be used for statistical inference.
\par In Theorem \ref{THEOREM 1}-\ref{THEOREM 3}, the choice of parameters $\rho_n$ and $\epsilon_n$ may not be optimal. In particular, the conditions guaranteeing asymptotic normality seem to be stronger than those guaranteeing collinearity. Furthermore, the asymptotic nature of the estimator is contingent not solely on the growth rate of the parameter but also on the configuration of the parameter. It would be interesting to investigate whether these boundary conditions could be relaxed.
\par We assume that the edges are binary edges following a Bernoulli distribution. In real-world networks, edges may not only take binary values but also be weighted. The theoretical analysis of weighted edges is worth further exploration. In the future, we will consider generalizing the analysis to other probability distributions, such as the Poisson distribution \cite{xu2023network}. Additionally, this paper makes the following assumptions: First, all parameters are bounded. Second, the generated network exhibits a dense structure, owing to the fact that the edge generation probability is significantly non-zero. However, real-world network data often show sparsity. Thus, how to analyze the asymptotic properties of differential privacy of network models incorporating covariates when data is sparse warrants further in-depth research and discussion . In future research, we will also consider generalizing the asymptotic theoretical approach for deriving differential privacy in this paper to a class of network models containing covariates to provide a unified theoretical framework.

\clearpage
\section*{Acknowledgements}
	Luo's research was partly supported by Humanities and Social Science Fund of Ministry of Education of China (Grant Number:24YJC910006)and by the Fundamental Research Funds for the Central Universities(South-Central MinZu University(Grant Number:CZQ24018)). And Qin's research was partly supported by National Natural Science Foundationon of China (Grant Number:12371261)

\end{document}